\documentclass[preprint, 1p]{elsarticle}

\usepackage{lineno,hyperref}
\modulolinenumbers[5]

\usepackage{amsmath}
\usepackage{amssymb}
\usepackage{cases}

\usepackage{mathptmx}      % use Times fonts if available on your TeX system

\usepackage{natbib}
\usepackage[english]{babel}
\usepackage{subfig}
\usepackage{graphicx}

\usepackage[active]{srcltx} % use TeX-souce-specials-mode

\usepackage[bbgreekl]{mathbbol}
\usepackage{bm} % sophisticated \boldsymbol
\usepackage{MnSymbol} % \lsem, \rsem, tensor product :
\usepackage{units}
\usepackage{tensor}
\usepackage{accents}

\usepackage{todonotes}

\DeclareMathOperator{\divergence}{div}

\DeclareMathOperator{\Tr}{Tr}

\newcommand{\exponential}[1]{\ensuremath{{\mathrm e}^{#1}}}

\newcommand{\reference}{\mathrm{ref}}

\newcommand{\bydefinition}{\mathrm{def}}
\newcommand{\traceless}[1]{{#1}_{\delta}}

\newcommand{\diff}{\mathrm{d}}
\newcommand{\Diff}[1][]{\mathrm{D}_{#1}} % For Frechet and Gateaux derivative
\renewcommand{\vec}[1]{\ensuremath{\mathbf{#1}}}

\makeatletter
\@ifpackageloaded{bm}% 
{\renewcommand{\vec}[1]{\ensuremath{\bm{#1}}}%
}{%
\relax% do nothing
}
\makeatother
\newcommand{\tensorq}[1]{\ensuremath{\mathbb{#1}}}      % tensorial quantity
\newcommand{\transpose}[1]{#1^\top}

\newcommand{\inverse}[1]{#1^{-1}}

\newcommand{\identity}{\ensuremath{\tensorq{I}}}

\newcommand{\cstress}{\tensorq{T}}

\newcommand{\fgrad}{\tensorq{F}}

\newcommand{\lcg}{\tensorq{B}}

\makeatletter
\@ifpackageloaded{bm}% 
{%
}{%

}

\@ifpackageloaded{bm}%
{%
 
}{%

}

\@ifpackageloaded{bm}%
{%
}{%

}
\makeatother

\newcommand{\generictensor}{{\tensorq{A}}}
\newcommand{\vecv}{\ensuremath{\vec{v}}}
\newcommand{\gradv}{\ensuremath{\nabla \vecv}}
\newcommand{\gradasym}{\ensuremath{\tensorq{W}}}
\newcommand{\gradsym}{\ensuremath{\tensorq{D}}}

\newcommand{\gradvl}{\ensuremath{\tensorq{L}}}

\newcommand{\R}{\ensuremath{{\mathbb R}}}
\newcommand{\ienergy}{\ensuremath{e}} % specific internal energy (energy per unit mass)
\newcommand{\fenergy}{\ensuremath{\psi}} % specific free energy
\newcommand{\entropy}{\ensuremath{\eta}} % specific entropy
\newcommand{\temp}{\ensuremath{\theta}} % temperature
\newcommand{\mns}{\ensuremath{m}} % mean normal stress
\newcommand{\nettenergy}{\ensuremath{E}_{\mathrm{tot}}} % net total energy
\newcommand{\netentropy}{\ensuremath{S}} % net entropy
\newcommand{\cheatvolref}{\ensuremath{c_{\mathrm{V}, \reference}}} % reference value

\newcommand{\efluxc}{\vec{j}_{e}} % energy flux, current configuration
\newcommand{\entfluxc}{\vec{j}_{\entropy}} % entropy flux, current configurtion 
\newcommand{\entprodc}{\xi} % entorpy production, current configuration
\newcommand{\entprodctemp}{\zeta} % entorpy production times temperature, current configuration

\newcommand{\pd}[2]{\ensuremath{\frac{\partial {#1}}{\partial {#2}}}}

\newcommand{\dd}[2]{\ensuremath{\frac{\diff {#1}}{\diff {#2}}}}

\newcommand{\ddd}[2]{\ensuremath{\frac{\diff^2 {#1}}{\diff {#2}^2}}}

\newcommand{\fid}[1]{\ensuremath{\accentset{\triangledown}{#1}}}

\newcommand{\gfid}[1]{\ensuremath{\accentset{\medsquare}{#1}}}

\makeatletter
\@ifpackageloaded{tensor}% tensor is a package for a better typesetting of tensors
{

}{%

}
\makeatother

\makeatletter
\@ifpackageloaded{tensor}% tensor is a package for a better typesetting of tensors
{

}{%

}
\makeatother

\newcommand{\norm}[2][]{\ensuremath{\left\|#2\right\|_{#1}}}
\newcommand{\absnorm}[1]{\ensuremath{\left|#1\right|}}

\makeatletter
\@ifundefined{volume}{%
}%
{%
}
\makeatother

\newcommand{\cvolumee}{\diff \mathrm{v}}

\makeatletter
\@ifpackageloaded{MnSymbol} % : as binary operator needs MnSymbol package
{
\newcommand{\tensordot}[2]{\ensuremath{#1 \vdotdot #2}} 
}{%
\newcommand{\tensordot}[2]{\ensuremath{#1 : #2}} 
}
\makeatother
\newcommand{\vectordot}[2]{\ensuremath{#1 \bullet #2}}

\newcommand{\sleb}[2]{\ensuremath{L}^{#1} \left(#2 \right)}             % Lebesgue space

\newcommand{\tensorf}[1]{{\mathfrak{#1}}}

\newcommand{\entropyrellp}{{\mathcal S}}

\newcommand{\energyrellp}{{\mathcal E}}

\newcommand{\tempbdr}{\temp_{\mathrm{bdr}}}

\newcommand{\tempref}{\ensuremath{\temp_{\reference}}}

\newcommand{\kapparef}{\ensuremath{\kappa_{\reference}}}

\newcommand{\fgradI}{\ensuremath{\fgrad_1}}

\newcommand{\fgradII}{\ensuremath{\fgrad_2}}
\newcommand{\lcgII}{\ensuremath{\lcg_{\mathrm{2}}}}

\newcommand{\lcgGSII}{\ensuremath{\lcg_{\mathrm{2}, \mathrm{GS}}}}

\numberwithin{equation}{section}
\let\cite\citet

\addto\captionsenglish{}

\begin{document}

\begin{frontmatter}

\title{Unconditional finite amplitude stability of a viscoelastic fluid in a mechanically isolated vessel with spatially non-uniform wall temperature}


\author[charles]{Mark Dostal\'{\i}k\fnref{myfootnote1}}
\ead{dostalik@karlin.mff.cuni.cz}
\fntext[myfootnote1]{Mark Dostal\'{\i}k has been supported by Charles University Research program No. UNCE/SCI/023 and Charles University Grant Agency, grant number~1652119.}

\author[charles]{V\'{\i}t Pr\r{u}\v{s}a\fnref{myfootnote2}\corref{mycorrespondingauthor}}
\fntext[myfootnote2]{V\'{\i}t Pr\r{u}\v{s}a thanks the Czech Science Foundation, grant number~18-12719S, for its support.}
\ead{prusv@karlin.mff.cuni.cz}

\author[heidelberg]{Judith Stein}
\ead{judith.stein@iwr.uni-heidelberg.de}

\address[charles]{Charles University, Faculty of Mathematics and Physics\\
Sokolovsk\'a 83, Praha, CZ 186 75, Czech Republic}
\address[heidelberg]{Universit\"at Heidelberg, Institute of Applied Mathematics\\
Im Neuenheimer Feld 205, Heidelberg, DE 69120, Germany}

\begin{abstract}
  We investigate finite amplitude stability of spatially inhomogeneous steady state of an incompressible viscoelastic fluid which occupies a mechanically isolated vessel with walls kept at spatially non-uniform temperature. For a wide class of incompressible viscoelastic models including the Oldroyd-B model, the Giesekus model, the FENE-P model, the Johnson--Segalman model, and the Phan--Thien--Tanner model we prove that the steady state is stable subject to any finite perturbation.

%%% Local Variables:
%%% mode: latex
%%% TeX-master: "../viscoelastic-stability-general"
%%% End:

\end{abstract}

\begin{keyword}
  finite amplitude stability\sep thermodynamically open system\sep non-equilibrium steady state\sep heat conducting fluid \sep viscoelastic fluid 
\MSC[2010] 
35Q35\sep  % PDEs in connection with fluid mechanics
35B35\sep  % Stability
37L15\sep  % Stability problems
76A10  % Viscoelastic fluids
% http://www.ams.org/msc/
\end{keyword}

\end{frontmatter}

\section{Introduction}
\label{sec:introduction}
We are interested in the long time behaviour of a fluid occupying a vessel that is mechanically isolated and that is allowed to exchange thermal energy with the surroundings. (The temperature boundary condition is an inhomogeneous Dirichlet boundary condition.) If no external forces are present, then one expects that the fluid in the vessel comes to the rest state as time goes to infinity. Moreover, the stability is expected to be unconditional, that is the rest state should be attained irrespective of the initial state of the fluid. The question is whether one can prove that such a long time behaviour is indeed implied by the corresponding governing equations.

Since the walls of the vessel are kept at a given~\emph{spatially nonuniform} temperature, the corresponding steady state is a~\emph{spatially inhomogeneous} solution to the governing equations, and the entropy is being produced (at a constant rate) at the steady state. Consequently, from the thermodynamic perspective the steady state is a~\emph{non-equilibrium} (entropy producing) steady state of a \emph{thermodynamically~open system}. This makes the analysis of the long time behaviour difficult as we cannot use methods developed for~thermodynamically isolated systems or for systems that are immersed in a thermal bath (spatially homogeneous temperature boundary condition), see~\cite{coleman.bd:on}, \cite{gurtin.me:thermodynamics*1,gurtin.me:thermodynamics} and later developments.

Recently, the issue of application of thermodynamically based methods in the stability a\-na\-ly\-sis of spatially inhomogeneous steady states has been discussed by~\cite{bulcek.m.malek.j.ea:thermodynamics}, where the authors have also proposed a systematic thermodynamically based approach to the stability problem. The approach proposed by~\cite{bulcek.m.malek.j.ea:thermodynamics} has been then used by~\cite{dostalk.m.prusa.v.rajagopal.kr:unconditional}, who have investigated the same stability problem as in the current contribution, but who have considered the Navier--Stokes--Fourier fluid (incompressible viscous heat conducting fluid).

Using minimal assumptions concerning the behaviour of the dissipative heating term in the evolution equation for temperature, \cite{dostalk.m.prusa.v.rajagopal.kr:unconditional} have shown that the corresponding spatially inhomogenoeus steady state is indeed unconditionally stable. In the present contribution we follow the approach by~\cite{dostalk.m.prusa.v.rajagopal.kr:unconditional}, and \emph{we generalise the findings by~\cite{dostalk.m.prusa.v.rajagopal.kr:unconditional} to include a variety of viscoelastic models}.

The analysis by~\cite{dostalk.m.prusa.v.rajagopal.kr:unconditional} has been based on two qualitative properties of the Navier--Stokes--Fourier model. First, the dissipative heating term in the evolution equation for the temperature must be a positive and integrable quantity. With a minimal effort we can show that this property is valid also for the considered viscoelastic rate-type models. Second, a norm of the velocity field must decay to zero at an exponential rate. This property is more complicated to show for the viscoelastic rate-type models, and its proof constitutes the main body of the current contribution. (In fact only show that the norm of the velocity field is bounded from above by an exponentially decaying function, but this is sufficient for the stability.) Once we show that the essential qualitative properties are preserved for viscoelastic rate-type fluids, it is straightforward to follow~\cite{dostalk.m.prusa.v.rajagopal.kr:unconditional}, and show the decay of the temperature perturbations.  

In particular, we show that stability of the spatially inhomogeneous non-equi\-librium steady state is indeed implied by the corresponding governing equations for the standard \hbox{Oldroyd-B} model, see~\cite{oldroyd.jg:on}, the Giesekus model, see~\cite{giesekus.h:simple}, the FENE-P model, see~\cite{bird.rb.dotson.pj.ea:polymer} and~\cite{keunings.r:on}, the Johnson--Segalman model, see~\cite{johnson.mw.segalman.d:model}, and the Phan--Thien--Tanner model, see~\cite{phan-thien.n.tanner.ri:new} and~\cite{phan-thien.n:non-linear}.

%\section{Preliminaries}
%\label{sec:preliminaries}
%\begin{lemma}[Decay of integrable functions]
%  \label{lm:1}
%  \begin{subequations}
%    \label{eq:5}
%  Let $y: [0, +\infty) \mapsto \R^+$ be a continuous non-negative function such that
%  \begin{equation}
%    \label{eq:6}
%    \int_{\tau=0}^{+\infty} y(\tau) \, \diff \tau \leq C_1,
%  \end{equation}
%  where $C_1$ is a constant. Moreover, let for all $s, t \in [0, + \infty)$, $t>s$,
%  \begin{equation}
%    \label{eq:7}
%    y(t) - y(s) \leq \int_{\tau=s}^t f(y(\tau)) \, \diff \tau +  \int_{\tau=s}^t h(\tau) \, \diff \tau,
%  \end{equation}
%  hold, where $f$ is a non-decreasing function from $\R^+$ to $\R^+$ and $h$ is a non-negative function such that
%  \begin{equation}
%    \label{eq:8}
%    \int_{\tau=0}^{+\infty} h(\tau) \, \diff \tau \leq C_2,
%  \end{equation}
%  where $C_2$ is a constant. Then
%  \begin{equation}
%    \label{eq:9}
%    \lim_{t \to + \infty} y(t) = 0.
%  \end{equation}
%  \end{subequations}
%\end{lemma}
%The lemma is taken from~\cite{zheng.s:asymptotic}, and its further generalisations can be found in~\cite{ramm.ag.hoang.ns:dynamical}. (See also {B}arb\u{a}lat lemma known in the optimal control theory, see~\cite{farkas.b.wegner.s:variations}.) Once we decide that Lemma~\ref{lm:1} or its suitable generalisation is the right tool for stability analysis, the only remaining task is to find the quantity $y(t)$ that satisfies the assumptions of the lemma, and that vanishes if and only if the temperature perturbation vanishes.

\section{General viscoelastic rate-type fluid}
\label{sec:general-viscoelastic-model}
Since the stability analysis will be based on thermodynamical concepts, we need to recall some facts regarding the themodynamic basis of the viscoelastic rate-type models for incompressible fluids. We present the derivation of a general thermodynamically consistent model which, among others, includes the Oldroyd-B model, the Giesekus model, the FENE-P model, the Johnson--Segalman model, and the Phan--Thien--Tanner model.

The derivation outlined below follows the procedure introduced by~\cite{rajagopal.kr.srinivasa.ar:thermodynamic}. The method is purely phenomenological and is based on the characterisation of the energy storage and entropy production mechanisms in the material. Specifically, we are interested in the identification of the specific Helmholtz free energy~$\fenergy$, see Section~\ref{sec:helmholtz-free-energy}, and the entropy production $\entprodc$, see Section~\ref{sec:entropy-production}. In the specific case of viscoelastic fluids we further apply a decomposition of its motion according to the dissipative and elastic response of the material. We \emph{virtually} split the deformation from the initial configuration to the current configuration into the deformation of the intermediate configuration, and to the instantaneous elastic deformation from the intermediate configuration to the current configuration, see Figure \ref{fig:viscoelastic-kinematics}. Such a decomposition of the total deformation to elastic and dissipative part then leads to certain kinematical identities that can be exploited in the derivation of the model.

\subsection{Kinematics}
\label{sec:kinematics}
\begin{figure}[h]
  \centering
  \includegraphics[width=0.5\textwidth]{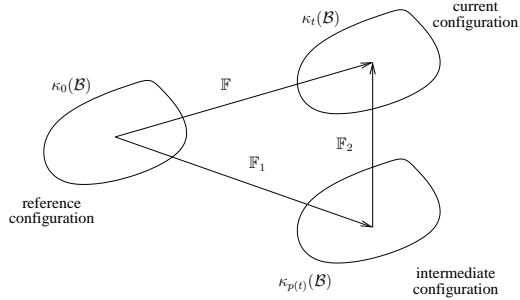}
  \caption{General decomposition of deformation gradient.}
  \label{fig:viscoelastic-kinematics}
\end{figure}

Let us concentrate on the decomposition of the motion of a viscoelastic body as depicted in Figure \ref{fig:viscoelastic-kinematics}. (For details see also~\cite{dostalk.m.prusa.v.ea:on}.) The total deformation gradient $\fgrad$ can be seen as a composition of two deformations
\begin{equation}
  \label{eq:fgrad-composition}
  \fgrad
  =
  \fgradII \fgradI,
\end{equation} 
where $\fgradI$ and $\fgradII$ are the deformation gradients of the partial deformations. Let us introduce the left Cauchy--Green tensor $\lcgII$ associated with the elastic response of the material via the relation
\begin{equation}
  \label{eq:lcgII}
  \lcgII
  =_{\bydefinition}
  \fgradII \transpose{\fgradII}.
\end{equation}
Tensor $\lcgII$ provides us a characterisation of the instantaneous elastic part of the deformation, and as we shall see in Section \ref{sec:constitutive-relations} it constitutes an additional ``elastic'' part of the Cauchy stress tensor.

The described decomposition yields viscoelastic models with the evolution equation containing the upper convected derivative
\begin{equation}
  \label{eq:upper-convected}
  \fid{\generictensor}
  =_{\bydefinition}
  \dd{\generictensor}{t} - \gradvl \generictensor - \generictensor \transpose{\gradvl},
\end{equation}
where
%\begin{equation}
%  \label{eq:material-derivative}
$
\dd{}{t} =_{\bydefinition} \pd{}{t} + \vectordot{\vecv}{\nabla},
$
% \end{equation}
denotes the material derivative, $\vecv$ denotes the spatial velocity, $\gradvl =_{\bydefinition} \gradv$ denotes the velocity gradient, and $\gradsym =_{\bydefinition} \frac{1}{2} \left( \gradvl + \transpose{\gradvl} \right)$ denotes the symmetric part of the velocity gradient. This setting is thus able to incorporate the standard Oldroyd-B model, the Giesekus model, and the FENE-P model.

However, the evolution equations for the Johnson--Segalman model and the Phan--Thien--Tanner model contain the so-called Gordon--Schowalter derivative
\begin{equation}
  \label{eq:gordon-schowalter}
  \gfid{\generictensor}
  =_{\bydefinition}
  \dd{\generictensor}{t} - a \left( \gradsym \generictensor + \generictensor \gradsym \right) - \left( \gradasym \generictensor + \generictensor \transpose{\gradasym} \right),
\end{equation}
where $a \in [-1, 1]$ and $\gradasym =_{\bydefinition} \frac{1}{2} \left( \gradvl - \transpose{\gradvl} \right)$ denotes the skew-symmetric part of the velocity gradient. Although the Gordon--Schowalter derivative is in general different from the upper convected derivative, it can be also obtained using the decomposition described above. However a generalisation of the decomposition~\eqref{eq:fgrad-composition} is needed. In principle one has to articulate the concept of ``non-affine'' motion introduced in \cite{johnson.mw.segalman.d:model}, see~\cite{dostalk.m.prusa.v.ea:on} for details. The generalised decomposition yields a different tensorial quantity associated with the additional ``elastic'' part of the Cauchy stress tensor. We denote this quantity by $\lcgGSII$. For a thorough analysis of the motion of a viscoelastic body in this generalised setting and interpretation of the tensorial quantity $\lcgGSII$, see \cite{dostalk.m.prusa.v.ea:on}.

Note that for $a = 1$ the Gordon--Schowalter derivative \eqref{eq:gordon-schowalter} reduces to the upper convected derivative and the tensorial quantity $\lcgGSII$ is simply recast to $\lcgII$. In the following, we shall thus be using the general notation $\lcgGSII$ for the additional tensorial quantity in the Cauchy stress tensor. For the models containing the upper convected derivative we then simply set $a = 1$ and use the notation $\lcgII$ instead of $\lcgGSII$.

\subsection{Helmholtz free energy}
\label{sec:helmholtz-free-energy}
We consider the specific Helmholtz free energy in the form
\begin{equation}
  \label{eq:fenergy}
  \fenergy 
  =_{\bydefinition}
  \fenergy_0(\temp)
  +
  \fenergy_1(\lcgGSII),
\end{equation}
where the thermal part $\fenergy_0$ is given by a simple formula (the symbols $\cheatvolref$ and $\tempref$ denote the specific heat capacity at constant volume and the reference temperature)
\begin{equation}
  \label{eq:fenergy-0}
  \fenergy_0
  =_{\bydefinition}
  -\cheatvolref \temp \left[ \ln \left( \frac{\temp}{\tempref} \right) - 1 \right],
\end{equation}
and $\fenergy_1$ satisfies the following set of requirements
\begin{subequations}
  \label{eq:fenergy-1-assumptions}
  \begin{gather}
    \label{eq:fenergy-1-assumption-zero-value}
    \fenergy_1 (\lcgGSII) \ge 0, \qquad \fenergy_1(\lcgGSII) = 0 \iff \lcgGSII = \identity,
    \\
    \label{eq:fenergy-1-assumption-derivative-zero-value}
    \pd{\fenergy_1}{\lcgGSII}(\lcgGSII) = \tensorq{0} \iff \lcgGSII = \identity,
    \\
    \label{eq:fenergy-1-assumption-commutativity}
    \lcgGSII \pd{\fenergy_1}{\lcgGSII}(\lcgGSII) = \pd{\fenergy_1}{\lcgGSII}(\lcgGSII) \lcgGSII.
  \end{gather}
\end{subequations}
%Note that assumption \eqref{eq:fenergy-1-assumption-commutativity} is equivalent to the product $\lcgGSII \pd{\fenergy_1}{\lcgGSII}(\lcgGSII)$ being a symmetric positive definite tensor.
%
(The commutative property~\eqref{eq:fenergy-1-assumption-commutativity} is immediately granted for the isotropic material.) The model-dependent quantity $\fenergy_1$ is specified in~\ref{sec:admissible-viscoelastic-models} for all the viscoelastic models mentioned in Section~\ref{sec:general-viscoelastic-model}. In the same section we also verify that the structural assumptions~\eqref{eq:fenergy-1-assumptions} are fulfilled for all considered models. Using the standard thermodynamic relations for the specific entropy $\entropy$  and the specific internal energy~$\ienergy$ 
\begin{subequations}
  \label{eq:thermodynamic-relations}
  \begin{align}
    \label{eq:thermodynamic-relations-entropy}
    \entropy &= - \pd{\fenergy}{\temp},
    \\
    \label{eq:thermodynamic-relations-energy}
    \ienergy &= \fenergy + \temp \entropy,
  \end{align}
\end{subequations}
together with the general evolution equation for the internal energy
\begin{equation}
  \label{eq:inergy-evolution-equation}
  \rho \dd{\ienergy}{t} = \tensordot{\cstress}{\gradsym} - \divergence \efluxc,
\end{equation}
we can derive an evolution equation for the specific entropy. (Here $\rho$ denotes density, $\cstress$ denotes the Cauchy stress tensor, and $\efluxc$ denotes the non-mechanical contribution to the energy flux.) Indeed, by taking the material derivative of \eqref{eq:thermodynamic-relations-energy} and exploiting the relations \eqref{eq:thermodynamic-relations-entropy} and \eqref{eq:inergy-evolution-equation} we arrive at
\begin{equation}
  \label{eq:1}
  \rho \dd{\entropy}{t} + \divergence \left( \frac{\efluxc}{\temp} \right)
  =
  \frac{1}{\temp}
  \left(
    \tensordot{\traceless{\cstress}}{\traceless{\gradsym}}
    -
    \rho \tensordot{\pd{\fenergy_1}{\lcgGSII}}{\dd{\lcgGSII}{t}}
    -
    \frac{\vectordot{\efluxc}{\nabla \temp}}{\temp}
  \right)
\end{equation}
Expressing the material derivative of $\lcgGSII$ via the formula for the Gordon--Schowalter derivative~\eqref{eq:gordon-schowalter}
\begin{equation}
  \label{eq:lcgGSII-gordon-schowalter}
  \dd{\lcgGSII}{t}
  =
  \gfid{\overline{\lcgGSII}} + a \left( \gradsym \lcgGSII + \lcgGSII \gradsym \right) + \left( \gradasym \lcgGSII + \lcgGSII \transpose{\gradasym} \right),
\end{equation}
and using the assumption \eqref{eq:fenergy-1-assumption-commutativity} we finally obtain
\begin{equation}
  \label{eq:entropy-evolution-equation-1}
  \rho \dd{\entropy}{t} + \divergence \left( \frac{\efluxc}{\temp} \right)
  =
  \frac{1}{\temp}
  \left\{
    \tensordot{\left[ \traceless{\cstress} - 2 \rho a \traceless{\left( \lcgGSII \pd{\fenergy_1}{\lcgGSII} \right)} \right]}{\traceless{\gradsym}}
    -
    \rho \tensordot{\pd{\fenergy_1}{\lcgGSII}}{\gfid{\overline{\lcgGSII}}}
    -
    \frac{\vectordot{\efluxc}{\nabla \temp}}{\temp}
  \right\}.
\end{equation}

\subsection{Entropy production}
\label{sec:entropy-production}
In order to identify the constitutive relations we want to ``compare'' equation \eqref{eq:entropy-evolution-equation-1} with the general evolution equation for entropy 
\begin{equation}
  \label{eq:entropy-evolution-equation-2}
  \rho \dd{\entropy}{t} + \divergence \entfluxc
  =
  \entprodc,
\end{equation}
where $\entfluxc$ denotes the entropy flux and the entropy production $\entprodc$ is given by
\begin{equation}
  \label{eq:entprod}
  \entprodc 
  =_{\bydefinition} 
  \frac{1}{\temp}
  \left(
    \entprodctemp_{\mathrm{th}}
    +
    \entprodctemp_{\mathrm{mech}}
  \right),
\end{equation}
where we have introduced the notation
\begin{subequations}
  \label{eq:zetas}
  \begin{align}
    \label{eq:zeta-thermal}
    \entprodctemp_{\mathrm{th}}
    &=_{\bydefinition}
    \kapparef \frac{\vectordot{\nabla \temp}{\nabla \temp}}{\temp},
    \\
    \label{eq:zeta-mechanical}
    \entprodctemp_{\mathrm{mech}}
    &=_{\bydefinition}
    2 \nu(\temp) \tensordot{\gradsym}{\gradsym}
    +
    \rho \frac{\mu}{\nu_1(\temp)} \tensordot{\pd{\fenergy_1}{\lcgGSII} (\lcgGSII)}{\tensorf{f}(\lcgGSII)}.
  \end{align}
\end{subequations}
Here, the symbol $\kapparef$ denotes the thermal conductivity, the material coefficient $\mu$ is a positive constant while the material coefficients $\nu$, $\nu_1$ are assumed to be positive functions of temperature. We require~$\nu$ to be bounded from below, and~$\nu_1$ to be bounded from above. Further, we assume that the tensorial function $\tensorf{f} : \R_{>}^{3 \times 3} \to \R_{>}^{3 \times 3}$, where~$\R_{>}^{3 \times 3}$ denotes the space of symmetric positive definite $3 \times 3$ matrices, satisfies
\begin{subequations}
  \label{eq:f-assumptions}
  \begin{gather}
    \label{eq:f-zero}
    \tensorf{f}(\lcgGSII) = \tensorq{0}
    \iff
    \lcgGSII = \identity,
    \\
    \label{eq:pd-fenergy-1-f-nonnegative}
    \tensordot{\pd{\fenergy_1}{\lcgGSII} (\lcgGSII)}{\tensorf{f}(\lcgGSII)}
    \ge
    0,
    \\
    \label{eq:fenergy-1-f-stability-inequality}
    \fenergy_1 (\lcgGSII)
    \le
    C_{\tensorf{f}}
    \tensordot{\pd{\fenergy_1}{\lcgGSII} (\lcgGSII)}{\tensorf{f} (\lcgGSII)},
  \end{gather}
\end{subequations}
where $C_{\tensorf{f}}$ is a positive constant dependent on the choice of $\tensorf{f}$. See~\ref{sec:admissible-viscoelastic-models} for specification of the tensorial function $\tensorf{f}$ for all the viscoelastic models mentioned in Section~\ref{sec:general-viscoelastic-model}. In the same section we also verify that the structural assumptions~\eqref{eq:f-assumptions} are fulfilled for all considered models.

\subsection{Constitutive relations}
\label{sec:constitutive-relations}
Comparison of the entropy production $\entprodc$ given by \eqref{eq:entprod} with the right-hand side of \eqref{eq:entropy-evolution-equation-1} yields the sought constitutive relations for the mechanical quantities $\cstress$ and $\lcgGSII$,
\begin{subequations}
  \label{eq:constitutive-relations}
  \begin{align}
    \label{eq:constitutive-relations-cstress}
    \traceless{\cstress}
    &=
    2 \nu(\temp) \traceless{\gradsym} 
    + 
    2 \rho a \traceless{\left( \lcgGSII \pd{\fenergy_1}{\lcgGSII} \right)},
    \\
    \label{eq:constitutive-relations-lcg}
    \nu_1(\temp) \gfid{\overline{\lcgGSII}}
    &=
    - \mu \tensorf{f}(\lcgGSII),
    \\
    \intertext{as well as for the energy/entropy fluxes $\efluxc$ and $\entfluxc$,}
    \label{eq:constitutive-relation-eflux}
    \efluxc
    &=
    -\kapparef \nabla \temp,
    \\
    \label{eq:constitutive-relation-entflux}
    \entfluxc
    &=
    - \frac{\kapparef \nabla \temp}{\temp}.
  \end{align}
\end{subequations}

\subsection{Evolution equation for temperature}
\label{sec:evolution-equation-for-temperature}
It remains to derive the evolution equation for temperature. Using the relation $\entropy = - \pd{\fenergy}{\temp} = - \dd{\fenergy_0}{\temp}$ we can rewrite the evolution equation for entropy \eqref{eq:entropy-evolution-equation-2} as
\begin{equation}
  \label{eq:0}
  \rho \dd{}{t} \left( - \dd{\fenergy_0}{\temp} \right) + \divergence \entfluxc
  =
  \entprodc.
\end{equation}
Using the special choice of $\fenergy_0$ given by \eqref{eq:fenergy-0}, the postulated entropy production \eqref{eq:entprod}, and the constitutive relation for the entropy flux \eqref{eq:constitutive-relation-entflux} in \eqref{eq:0} then yields the evolution equation for temperature
\begin{equation}
  \label{eq:evolution-equation-for-temperature}
    \rho \cheatvolref \dd{\temp}{t}
    =
    \divergence(\kapparef \nabla \temp)
    +
    \entprodctemp_{\mathrm{mech}}.
  \end{equation}

  We note that the structure of the temperature evolution equation is the same both for the Navier--Stokes--Fourier fluid and for our general viscoelastic rate-type fluid. The two fluid models differ in the specification of the entropy production term $\entprodctemp_{\mathrm{mech}}$, see also~\ref{eq:zeta-mechanical}. \emph{Since the stability analysis done by~\cite{dostalk.m.prusa.v.rajagopal.kr:unconditional} required that the entropy production term $\entprodctemp_{\mathrm{mech}}$ is nonnegative and integrable in time and space, we see that this assumption is very likely to hold also for our general viscoelastic rate-type model. Consequently, one can conjecture that it would be possible to reuse much of the results obtained in~\cite{dostalk.m.prusa.v.rajagopal.kr:unconditional}.} As we shall see later, this is indeed the case. 

\section{Problem formulation}
\label{eq:problem-formulation}

\subsection{Governing equations and boundary conditions}
\label{sec:governing-equations}
Appealing to the derived constitutive relations~\eqref{eq:constitutive-relations} and the evolution equation for temperature \eqref{eq:evolution-equation-for-temperature} we see that the complete system of evolution equations describing the behaviour of our general viscoelastic rate-type fluid reads
\begin{subequations}
  \label{eq:governing-equations}
  \begin{align}
    \label{eq:incompressibility-condition}
    \divergence \vecv 
    &= 
    0,
    \\
    \label{eq:linear-momentum-equation}
    \rho \dd{\vecv}{t}
    &=
    \nabla \mns + \divergence \left[ 2 \nu(\temp) \gradsym + 2 \rho a \traceless{\left( \lcgGSII \pd{\fenergy_1}{\lcgGSII} \right)} \right],
    \\
    \label{eq:lcg-equation} 
    \nu_1(\temp) \gfid{\overline{\lcgGSII}}
    &=
    - \mu \tensorf{f}(\lcgGSII),
    \\
    \label{eq:heat-equation}
    \rho \cheatvolref \dd{\temp}{t}
    &=
    \divergence(\kapparef \nabla \temp)
    +
    \entprodctemp_{\mathrm{mech}},
  \end{align}
\end{subequations}
where $\mns =_{\bydefinition} \frac{1}{3} \Tr \cstress$ denotes the mean normal stress.
%
%\subsection{Boundary conditions}
%\label{sec:boundary-conditions}
The evolution equations~\ref{eq:governing-equations} for the quadruple $\vec{W} = _{\bydefinition} \left[ \mns, \vecv, \lcgGSII, \temp \right]$ must be solved in the domain $\Omega$ that represents the closed vessel, while the boundary conditions on the vessel walls are\begin{subequations}
  \label{eq:boundary-conditions}
  \begin{align}
    \label{eq:boundary-conditions-velocity}
    \left. \vecv \right|_{\partial \Omega} &= \vec{0},
    \\
    \label{eq:boundary-conditions-temperature}
    \left. \temp \right|_{\partial \Omega} &= \tempbdr.
  \end{align}
\end{subequations}
The quantity $\tempbdr$ is a given nontrivial function of position.

\subsection{Problem of stability of the steady state}
\label{sec:problem-of-stability-of-the-steady-state}
The objective is to show that the perturbations $\widetilde{\vec{W}} =_{\bydefinition} \left[\widetilde{\mns}, \widetilde{\vecv}, \widetilde{\lcgGSII}, \widetilde{\temp} \right]$ to the steady state $\widehat{\vec{W}} =_{\bydefinition} \left[\widehat{\mns}, \widehat{\vecv}, \widehat{\lcgGSII}, \widehat{\temp} \right]$ vanish as time goes to infinity, that is
\begin{equation}
  \label{eq:perturbation-vanishes}
  \lim_{t \to +\infty} \widetilde{\vec{W}} = \vec{0}, 
\end{equation}
while the evolution of the quadruple $\vec{W} = \widehat{\vec{W}} + \widetilde{\vec{W}}$ is governed by evolution equations~\eqref{eq:governing-equations}.

\subsection{Spatially inhomogeneous non-equilibrium steady state}
\label{sec:spatially-inhomogeneous-non-equilibrium-steady-state}
In the non-equilibrium \emph{steady} state $\widehat{\vec{W}} =_{\bydefinition} [\widehat{\mns}, \widehat{\vecv}, \widehat{\lcgGSII}, \widehat{\temp}]$ the fluid is at rest $\widehat{\vecv} = \vec{0}$, and the tensorial quantity $\lcgGSII$ reduces to identity, that is $\widehat{\lcgGSII} = \identity$. This observation follows from \eqref{eq:lcg-equation} and the structural assumption~\eqref{eq:f-zero}. Further from \eqref{eq:linear-momentum-equation} and the assumption \eqref{eq:fenergy-1-assumption-derivative-zero-value} we obtain $\nabla \widehat{\mns} = 0$. Lastly, the temperature evolution equation \eqref{eq:heat-equation} implies that the steady temperature field $\widehat{\temp}$ solves
\begin{subequations}
  \begin{align}
    \label{eq:steady-state-temperature}
    0
    &=
    \divergence \left(\kapparef \nabla \widehat{\temp} \right),
    \\
    \label{eq:steady-state-boundary-conditions-temperature}
    \left. \temp \right|_{\partial \Omega} 
    &= 
    \tempbdr.
  \end{align}
\end{subequations}
The temperature field is thus given by the steady heat equation \eqref{eq:steady-state-temperature} with Dirichlet boundary condition \eqref{eq:steady-state-boundary-conditions-temperature}. If $\tempbdr$ is a nontrivial function of position, then $\widehat{\temp}$ is a spatially inhomogeneous bounded function. 

\subsection{Evolution equations for perturbations to the mechanical quantities}
\label{sec:evolution-equations-for-perturbations-to-the-mechanical-quantities}
Using the governing equations~\eqref{eq:governing-equations} it is straightforward to derive evolution equations for the perturbations $\widetilde{\vec{W}} =_{\bydefinition} [\widetilde{\mns}, \widetilde{\vecv}, \widetilde{\lcgGSII}, \widetilde{\temp}]$ to the steady state. The evolution equations for the mechanical quantities $\widetilde{\vec{v}}$ and $\widetilde{\lcgGSII}$ read
\begin{subequations}
  \label{eq:evolution-equations-perturbations}
  \begin{align}
    \label{eq:evolution-equations-perturbations-velocity}
    \rho \pd{\widetilde{\vecv}}{t}
    ={}&
    - 
    \rho \left( \vectordot{\widetilde{\vecv}}{\nabla} \right) \widetilde{\vecv}
    +
    \nabla \widetilde{\mns}
    +
    \divergence 
    \left[ 
      2 \nu(\widehat{\temp} + \widetilde{\temp}) \widetilde{\gradsym} + 2 \rho a \traceless{\left( (\identity + \lcgGSII) \pd{\fenergy_1}{\lcgGSII}(\identity + \lcgGSII) \right)}
    \right],
    \\
    \label{eq:evolution-equations-perturbations-lcg}
    \begin{split}
      \pd{\widetilde{\lcgGSII}}{t}
      ={}&
      - 
      \left( \vectordot{\widetilde{\vecv}}{\nabla} \right) \widetilde{\lcgGSII}
      +
      a 
      \left( 
        \widetilde{\gradsym} \widetilde{\lcgGSII} + \widetilde{\lcgGSII} \widetilde{\gradsym} 
      \right)
      +
      \widetilde{\gradasym} \widetilde{\lcgGSII}
      +
      \widetilde{\lcgGSII} \transpose{\widetilde{\gradasym}}
      +
      2 a \widetilde{\gradsym}
      \\
      &-
      \frac{\mu}{\nu_1(\widehat{\temp} + \widetilde{\temp})} \tensorf{f}(\identity + \widetilde{\lcgGSII}).
    \end{split}
  \end{align}
(In the derivation of~\eqref{eq:evolution-equations-perturbations-velocity} we have exploited the assumption \eqref{eq:fenergy-1-assumption-derivative-zero-value}.) Furthermore, the evolution equation for the temperature perturbation $\widetilde{\temp}$ reads
\begin{equation}
  \label{eq:3}
  \rho \cheatvolref \pd{\widetilde{\temp}}{t}
  +
  \rho \cheatvolref \vectordot{\widetilde{\vec{v}}}{\left[ \nabla \left(\widehat{\temp} + \widetilde{\temp}\right) \right]}
  =
  \divergence
  \left(
    \kapparef
    \nabla
    \widetilde{\temp}
  \right)
  +
  \entprodctemp_{\mathrm{mech}}\left(\widehat{\vec{W}} + \widetilde{\vec{W}}\right).
\end{equation}
\end{subequations}

\section{Thermodynamically motivated construction of a Lyapunov type functional}
\label{eq:thermodynamically-motivated-functional}
The stability is investigated using the concepts introduced in~\cite{bulcek.m.malek.j.ea:thermodynamics} and~\cite{dostalk.m.prusa.v.rajagopal.kr:unconditional}.

\subsection{Construction of the functional}
\label{sec:construction-of-the-functional}
Following~\cite{bulcek.m.malek.j.ea:thermodynamics} we define Lyapunov type functional $\mathcal{V}_{\mathrm{neq}}$ as
\begin{equation}
  \label{eq:functional-definition}
  \mathcal{V}_{\mathrm{neq}}
  \left(
    \left.
      \widetilde{\vec{W}}
    \right\|
    \widehat{\vec{W}}
  \right)
  =_{\bydefinition}
  -
  \left[
    \entropyrellp_{\widehat{\temp}}(\left. \widetilde{\vec{W}} \right\| \widehat{\vec{W}})
    -
    \energyrellp  (\left. \widetilde{\vec{W}} \right\| \widehat{\vec{W}})
  \right]
  ,
\end{equation}
where
\begin{subequations}
  \label{eq:rells}
  \begin{align}
    \label{eq:entropyrellp-definition}
    \entropyrellp_{\widehat{\temp}}(\left. \widetilde{\vec{W}} \right\| \widehat{\vec{W}})
    &=_{\bydefinition}
    \netentropy_{\widehat{\temp}}(\widehat{\vec{W}} + \widetilde{\vec{W}})
    -
    \netentropy_{\widehat{\temp}}(\widehat{\vec{W}})
    -
    \Diff \netentropy_{\widehat{\temp}} (\widehat{\vec{W}})
    \left[
      \widetilde{\vec{W}}
    \right],
    \\
    \label{eq:energyrellp-definition}
    \energyrellp  (\left. \widetilde{\vec{W}} \right\| \widehat{\vec{W}})
    &=_{\bydefinition}
    \nettenergy(\widehat{\vec{W}} + \widetilde{\vec{W}})
    -
    \nettenergy(\widehat{\vec{W}})
    -
    \Diff \nettenergy (\widehat{\vec{W}})
    \left[
      \widetilde{\vec{W}}
    \right],
  \end{align}
\end{subequations}
and the rescaled net entropy $\netentropy_{\widehat{\temp}}$ and the net total energy $\nettenergy$ are given by the formulae
\begin{subequations}
  \label{eq:functional-components}
  \begin{align}
    \label{eq:netentropy}
    \netentropy_{\widehat{\temp}}(\vec{W})
    &=_{\bydefinition}
    \int_{\Omega} \! \rho \widehat{\temp} \entropy (\vec{W}) \, \cvolumee
    =
    - 
    \int_{\Omega} \! 
      \rho \widehat{\temp} \dd{\fenergy_0}{\temp} (\temp)
    \cvolumee,
    \\
    \label{eq:netenergy}
    \nettenergy(\vec{W})
    &=_{\bydefinition}
    \int_{\Omega} \! 
      \left[
        \rho \ienergy(\vec{W}) + \frac{1}{2} \rho \absnorm{\vecv}^2 
      \right]
    \cvolumee
    =
    \int_{\Omega} \! 
      \rho
      \left[
        \fenergy_0(\temp)
        + 
        \fenergy_1(\lcgGSII)
        -
        \temp \dd{\fenergy_0}{\temp} (\temp)
        +
        \frac{1}{2} \rho \absnorm{\vecv}^2
      \right]
    \cvolumee,
  \end{align}
\end{subequations}
where $\ienergy$ denotes the specific internal energy, $\entropy$ denotes the specific entropy, and where have exploited thermodynamic relations \eqref{eq:thermodynamic-relations}. In \eqref{eq:rells}, the symbols $\Diff \netentropy_{\widehat{\temp}} (\widehat{\vec{W}}) [\widetilde{\vec{W}}]$ and $\Diff \nettenergy (\widehat{\vec{W}}) [\widetilde{\vec{W}}]$ denote the G\^{a}teaux derivative of the given functional at point $\widehat{\vec{W}}$ in the direction $\widetilde{\vec{W}}$. It particular, we have
\begin{subequations}
  \label{eq:gateaux-derivatives}
  \begin{align}
    \label{eq:gateaux-derivatives-netentropy}
    \Diff \netentropy_{\widehat{\temp}} (\widehat{\vec{W}}) [\widetilde{\vec{W}}]
    &=
    \int_{\Omega} \! 
      \rho \widehat{\temp} \widetilde{\temp} \ddd{\fenergy_0}{\temp} (\widehat{\temp})
    \cvolumee,
    \\
    \Diff \nettenergy (\widehat{\vec{W}}) [\widetilde{\vec{W}}]
    &=
    \int_{\Omega}
      \rho
      \left[
        \tensordot{\widetilde{\lcgGSII}}{\pd{\fenergy_1}{\lcgGSII}}
        -
        \widehat{\temp} \widetilde{\temp} \ddd{\fenergy_0}{\temp} (\widehat{\temp})
        +
        \rho \vectordot{\widehat{\vecv}}{\widetilde{\vecv}}
      \right]
    \cvolumee,
  \end{align}
\end{subequations}

Consequently, it is straightforward to see that the formulae for the functionals $\entropyrellp_{\widehat{\temp}}$ and $\energyrellp$ read
\begin{subequations}
  \label{eq:functionals-formulae}
  \begin{align}
  \label{eq:entropyrellp}
  \entropyrellp_{\widehat{\temp}}(\left. \widetilde{\vec{W}} \right\| \widehat{\vec{W}})
  ={}&
  -
  \int_{\Omega}
    \rho \widehat{\temp}
    \left[
      \dd{\fenergy_0}{\temp} (\widehat{\temp} + \widetilde{\temp})
      -
      \dd{\fenergy_0}{\temp} (\widehat{\temp})
      -
      \widetilde{\temp} \ddd{\fenergy_0}{\temp} (\widehat{\temp})
    \right]
  \cvolumee,
  \\
  \label{eq:energyrellp}
  \begin{split}
    \energyrellp  (\left. \widetilde{\vec{W}} \right\| \widehat{\vec{W}})
    ={}&
    \int_{\Omega}
      \rho
      \bigg[
        \fenergy_0 (\widehat{\temp} + \widetilde{\temp})
        -
        \fenergy_0 (\widehat{\temp})
        -
        (\widehat{\temp} + \widetilde{\temp}) \dd{\fenergy_0}{\temp} (\widehat{\temp} + \widetilde{\temp})
        +
        \widehat{\temp} \dd{\fenergy_0}{\temp} (\widehat{\temp})
        +
        \widehat{\temp} \widetilde{\temp} \ddd{\fenergy_0}{\temp} (\widehat{\temp})
        \\
        &{\quad}+
        \fenergy_1 (\widehat{\lcgGSII} + \widetilde{\lcgGSII})
        -
        \fenergy_1 (\widehat{\lcgGSII})
        -
        \tensordot{\widetilde{\lcgGSII}}{\pd{\fenergy_1}{\lcgGSII} (\widehat{\lcgGSII})}
        +
        \frac{1}{2} \rho \absnorm{\widetilde{\vecv}}^2
      \bigg]
    \cvolumee,
  \end{split}
  \end{align}
\end{subequations}
hence the explicit formula for the functional $\mathcal{V}_{\mathrm{neq}}$ introduced in~\eqref{eq:functional-definition} reads
\begin{multline}
  \label{eq:functional-explicit-general}
  \mathcal{V}_{\mathrm{neq}}
  \left(
    \left.
      \widetilde{\vec{W}}
    \right\|
    \widehat{\vec{W}}
  \right)
  =
  \int_{\Omega}
    \rho
    \left[
      \fenergy_0 (\widehat{\temp} + \widetilde{\temp})
      -
      \fenergy_0 (\widehat{\temp})
      -
      \widetilde{\temp} \dd{\fenergy_0}{\temp} (\widehat{\temp} + \widetilde{\temp})
    \right]
  \cvolumee
  \\
  +
  \int_{\Omega}
    \rho
    \left[
      \fenergy_1 (\widehat{\lcgGSII} + \widetilde{\lcgGSII})
      -
      \fenergy_1 (\widehat{\lcgGSII})
      -
      \tensordot{\widetilde{\lcgGSII}}{\pd{\fenergy_1}{\lcgGSII} (\widehat{\lcgGSII})}
    \right]
  \cvolumee
  +
  \int_{\Omega}
    \frac{1}{2} \rho \absnorm{\widetilde{\vecv}}^2
  \cvolumee.
\end{multline}

For the subsequent stability analysis it is convenient to split the functional $\mathcal{V}_{\mathrm{neq}}$ into two parts
\begin{subequations}
  \label{eq:functional-parts-general}
  \begin{align}
    \label{eq:functional-thermal-part-general}
    \mathcal{V}_{\mathrm{th}}
    \left(
      \left.
        \widetilde{\vec{W}}
      \right\|
      \widehat{\vec{W}}
    \right)
    =_{\bydefinition}{}&
    \int_{\Omega}
      \rho
      \left[
        \fenergy_0 (\widehat{\temp} + \widetilde{\temp})
        -
        \fenergy_0 (\widehat{\temp})
        -
        \widetilde{\temp} \dd{\fenergy_0}{\temp} (\widehat{\temp} + \widetilde{\temp})
      \right]
    \cvolumee,
    \\
    \label{eq:functional-mechanical-part-general}
    \begin{split}
    \mathcal{V}_{\mathrm{mech}}
    \left(
      \left.
        \widetilde{\vec{W}}
      \right\|
      \widehat{\vec{W}}
    \right)
    =_{\bydefinition}{}&
    \int_{\Omega}
      \rho
      \left[
        \fenergy_1 (\widehat{\lcgGSII} + \widetilde{\lcgGSII})
        -
        \fenergy_1 (\widehat{\lcgGSII})
        -
        \tensordot{\widetilde{\lcgGSII}}{\pd{\fenergy_1}{\lcgGSII} (\widehat{\lcgGSII})}
      \right]
    \cvolumee
    \\
    &+
    \int_{\Omega}
      \frac{1}{2} \rho \absnorm{\widetilde{\vecv}}^2
    \cvolumee,
    \end{split}
  \end{align}
\end{subequations}
where $\mathcal{V}_{\mathrm{th}}$ shall be used to deal with the temperature perturbations $\widetilde{\temp}$, while $\mathcal{V}_{\mathrm{mech}}$ shall be used to deal with the perturbations to the mechanical quantities $\widetilde{\vecv}$ and $\widetilde{\lcgGSII}$. Note that in general $\mathcal{V}_{\mathrm{th}} \neq \entropyrellp_{\widehat{\temp}}$, $\mathcal{V}_{\mathrm{mech}} \neq \energyrellp$. However, if $\fenergy_0$ is chosen as in~\eqref{eq:fenergy-0}, then the corresponding functionals coincide.

Recall that so far we have considered the specific free energy in the general form \eqref{eq:fenergy}. However, in our specific case, $\fenergy_0$ is given by \eqref{eq:fenergy-0} and, moreover, in the steady non-equilibrium state we have $\widehat{\lcgGSII} = \identity$, which together with the assumptions \eqref{eq:fenergy-1-assumption-zero-value} and \eqref{eq:fenergy-1-assumption-derivative-zero-value} yields the final formula for $\mathcal{V}_{\mathrm{neq}}$
\begin{equation}
  \label{eq:functional-explicit-specific}
  \mathcal{V}_{\mathrm{neq}}
  \left(
    \left.
      \widetilde{\vec{W}}
    \right\|
    \widehat{\vec{W}}
  \right)
  =
  \int_{\Omega}
    \rho \cheatvolref \widehat{\temp}
    \left[
      \frac{\widetilde{\temp}}{\widehat{\temp}}
      -
      \ln \left( 1 + \frac{\widetilde{\temp}}{\widehat{\temp}} \right)
    \right]
  \cvolumee
  +
  \int_{\Omega}
    \rho \fenergy_1 (\identity + \widetilde{\lcgGSII})
  \cvolumee
  +
  \int_{\Omega}
    \frac{1}{2} \rho \absnorm{\widetilde{\vecv}}^2
  \cvolumee,
\end{equation}
along with
\begin{subequations}
  \label{eq:functional-parts}
  \begin{align}
    \label{eq:functional-thermal-part}
    \mathcal{V}_{\mathrm{th}}
    \left(
      \left.
        \widetilde{\vec{W}}
      \right\|
      \widehat{\vec{W}}
    \right)
    &=
    \int_{\Omega}
      \rho \cheatvolref \widehat{\temp}
      \left[
        \frac{\widetilde{\temp}}{\widehat{\temp}}
        -
        \ln \left( 1 + \frac{\widetilde{\temp}}{\widehat{\temp}} \right)
      \right]
    \cvolumee,
    \\
    \label{eq:functional-mechanical-part}
    \mathcal{V}_{\mathrm{mech}}
    \left(
      \left.
        \widetilde{\vec{W}}
      \right\|
      \widehat{\vec{W}}
    \right)
    &=
    \int_{\Omega}
      \rho \fenergy_1 (\identity + \widetilde{\lcgGSII})
    \cvolumee
    +
    \int_{\Omega}
      \frac{1}{2} \rho \absnorm{\widetilde{\vecv}}^2
    \cvolumee.
  \end{align}
\end{subequations}
It is straightforward to show that the functionals $\mathcal{V}_{\mathrm{neq}}$, $\mathcal{V}_{\mathrm{th}}$, and $\mathcal{V}_{\mathrm{mech}}$ are nonnegative and vanish if and only if the perturbation vanishes.

\subsection{Time derivative of the functional}
\label{sec:time-derivative-of-the-functional}
The time derivative of the thermal part $\mathcal{V}_{\mathrm{th}}$ of the constructed functional $\mathcal{V}_{\mathrm{neq}}$ has been already dealt with in~\cite{dostalk.m.prusa.v.rajagopal.kr:unconditional}, see Appendix~A therein, hence we will not repeat the lengthy algebraic manipulation here. (Note that although~\cite{dostalk.m.prusa.v.rajagopal.kr:unconditional} have considered the Navier--Stokes--Fourier fluid, their results regarding the thermal part of the proposed functional are applicable to viscoelastic rate-type fluids as well. This follows from the fact that the particular choice of the formula for the mechanical dissipation $\entprodctemp_{\mathrm{mech}}(\widehat{\vec{W}} + \widetilde{\vec{W}})$ has been inconsequential in the analysis by~\cite{dostalk.m.prusa.v.rajagopal.kr:unconditional}. See also Section~\ref{sec:evolution-equation-for-temperature} for a thorough discussion thereof.) The time derivative of $\mathcal{V}_{\mathrm{th}}$ is given by
\begin{multline}
  \label{eq:time-derivative-thermal-functional}
  \dd{\mathcal{V}_{\mathrm{th}}}{t}
  \left(
    \left.
      \widetilde{\vec{W}}
    \right\|
    \widehat{\vec{W}}
  \right)
  =
  -
  \int_{\Omega}
    \kapparef \widehat{\temp} \vectordot{\nabla \ln \left( 1 + \frac{\widetilde{\temp}}{\widehat{\temp}} \right)}{\nabla \ln \left( 1 + \frac{\widetilde{\temp}}{\widehat{\temp}} \right)}
  \, \cvolumee
  -
  \int_{\Omega}
    \rho \cheatvolref \ln \left( 1 + \frac{\widetilde{\temp}}{\widehat{\temp}} \right)
    \left(
      \vectordot{\widetilde{\vecv}}{\nabla \widehat{\temp}}  
    \right)  
  \, \cvolumee
  \\
  +
  \int_{\Omega}
    \frac{\widetilde{\temp}}{\widehat{\temp} + \widetilde{\temp}} \,
    \entprodctemp_{\mathrm{mech}}(\widehat{\vec{W}} + \widetilde{\vec{W}})
  \, \cvolumee.
\end{multline}

The formula for the time derivative of the mechanical part $\mathcal{V}_{\mathrm{mech}}$ of the constructed functional follows from the following manipulation. Direct differentiation under the integral sign yields
\begin{equation}
  \label{eq:time-derivative-mechanical-functional-proto}
  \dd{\mathcal{V}_{\mathrm{mech}}}{t}
  \left(
    \left.
      \widetilde{\vec{W}}
    \right\|
    \widehat{\vec{W}}
  \right)
  =
  \int_{\Omega}
    \rho 
    \tensordot{
      \pd{\fenergy_1 (\identity + \widetilde{\lcgGSII})}{\lcgGSII}
    }{
      \pd{\lcgGSII}{t}
    }
  \cvolumee
  +
  \int_{\Omega}
    \rho \vectordot{\widetilde{\vecv}}{\pd{\widetilde{\vecv}}{t}}
  \cvolumee.
\end{equation}
Using the evolution equation for the perturbation of left Cauchy--Green tensor \eqref{eq:evolution-equations-perturbations-lcg} the first term of \eqref{eq:time-derivative-mechanical-functional-proto} translates to
\begin{multline}
  \label{eq:term-1}
  \int_{\Omega}
    \rho 
    \tensordot{
      \pd{\fenergy_1 (\identity + \widetilde{\lcgGSII})}{\lcgGSII}
    }{
      \pd{\lcgGSII}{t}
    }
  \cvolumee
  =
  \int_{\Omega}
    2 \rho a 
    \tensordot{\left( (\identity + \lcgGSII) \pd{\fenergy_1}{\lcgGSII}(\identity + \lcgGSII) \right)}{\gradsym}
  \, \cvolumee
  \\
  -
  \int_{\Omega}
    \rho \frac{\mu}{\nu_1 (\widehat{\temp} + \widetilde{\temp})}
    \tensordot{
      \pd{\fenergy_1}{\lcgGSII} (\identity + \widetilde{\lcgGSII})
    }{
      \tensorf{f} (\identity + \widetilde{\lcgGSII})
    }
  \, \cvolumee,
\end{multline}
where we have used the assumption \eqref{eq:fenergy-1-assumption-commutativity} and the identity
\begin{equation}
  \label{eq:2}
  \int_{\Omega}
    \rho
    \tensordot{
      \pd{\fenergy_1 (\identity + \widetilde{\lcgGSII})}{\lcgGSII}
    }{
      \left(\vectordot{\widetilde{\vecv}}{\nabla}\right) \widetilde{\lcgGSII}
    }
  \, \cvolumee
  =
  \int_{\Omega}
    \vectordot{\widetilde{\vecv}}{\nabla \fenergy_1 \left( \identity + \lcgGSII \right)}
  \, \cvolumee
  =
  0.
\end{equation}
(The last equality follows from the Stokes theorem and from the fact that $\widetilde{\vecv}$ vanishes on the boundary.) Similarly, using the evolution equation for the velocity perturbation \eqref{eq:evolution-equations-perturbations-velocity}, the second term of \eqref{eq:time-derivative-mechanical-functional-proto} is recast into
\begin{equation}
  \label{eq:term-2}
  \int_{\Omega}
    \rho \vectordot{\widetilde{\vecv}}{\pd{\widetilde{\vecv}}{t}}
  \cvolumee
  =
  -
  \int_{\Omega}
    2 \nu(\widehat{\temp} + \widetilde{\temp}) \tensordot{\widetilde{\gradsym}}{\widetilde{\gradsym}}
  \, \cvolumee
  -
  \int_{\Omega}
    2 \rho a 
    \tensordot{\left( (\identity + \lcgGSII) \pd{\fenergy_1}{\lcgGSII}(\identity + \lcgGSII) \right)}{\gradsym}
  \, \cvolumee.
\end{equation}
Combining \eqref{eq:term-1} and \eqref{eq:term-2} in \eqref{eq:time-derivative-mechanical-functional-proto} and using the definition of $\entprodctemp_{\mathrm{mech}}$, see~\eqref{eq:zeta-mechanical}, we arrive at the final formula for the time derivative of $\mathcal{V}_{\mathrm{mech}}$,
\begin{equation}
  \label{eq:time-derivative-mechanical-functional}
  \dd{\mathcal{V}_{\mathrm{mech}}}{t}
  \left(
    \left.
      \widetilde{\vec{W}}
    \right\|
    \widehat{\vec{W}}
  \right)
  =
  -
  \int_{\Omega}
    \entprodctemp_{\mathrm{mech}} (\widehat{\vec{W}} + \widetilde{\vec{W}})
  \, \cvolumee.
\end{equation}

Consequently, equations \eqref{eq:time-derivative-thermal-functional} and \eqref{eq:time-derivative-mechanical-functional} yield the time derivative of the full functional $\mathcal{V}_{\mathrm{neq}}$
\begin{multline}
  \label{eq:time-derivative-functional}
  \dd{\mathcal{V}_{\mathrm{neq}}}{t}
  \left(
    \left.
      \widetilde{\vec{W}}
    \right\|
    \widehat{\vec{W}}
  \right)
  =
  -
  \int_{\Omega}
    \kapparef \widehat{\temp} \vectordot{\nabla \ln \left( 1 + \frac{\widetilde{\temp}}{\widehat{\temp}} \right)}{\nabla \ln \left( 1 + \frac{\widetilde{\temp}}{\widehat{\temp}} \right)}
  \, \cvolumee
  -
  \int_{\Omega}
    \rho \cheatvolref \ln \left( 1 + \frac{\widetilde{\temp}}{\widehat{\temp}} \right)
    \left(
      \vectordot{\widetilde{\vecv}}{\nabla \widehat{\temp}}  
    \right) 
  \, \cvolumee
  \\
  -
  \int_{\Omega}
    \frac{\widehat{\temp}}{\widehat{\temp} + \widetilde{\temp}} \,
    \entprodctemp_{\mathrm{mech}}(\widehat{\vec{W}} + \widetilde{\vec{W}})
  \, \cvolumee.
\end{multline}
In virtue of assumption \eqref{eq:pd-fenergy-1-f-nonnegative}, we know that the last term on the right-hand side of \eqref{eq:time-derivative-functional} is nonnegative. The only term whose sign is not known \emph{a priori} is the second term on the right-hand side of \eqref{eq:time-derivative-functional}. Its presence prohibits one from showing that the time derivative of the functional~$\mathcal{V}_{\mathrm{neq}}$ is, for a non-constant $\widehat{\temp}$, a nonpositive quantity. Consequently, $\mathcal{V}_{\mathrm{neq}}$ cannot directly serve as a genuine Lyapunov functional.

\subsection{Family of functionals $\mathcal{V}_{\mathrm{th}}^{m}$}
\label{sec:family-of-functionals-th}
As it has been shown in~\cite{dostalk.m.prusa.v.rajagopal.kr:unconditional}, the functional $\mathcal{V}_{\mathrm{th}}$ is insufficient to yield the asymptotic stability of the steady temperature field~$\widehat{\temp}$ via the Lyapunov method. This---rather technical---difficulty can be dealt with by introducing a new temperature scale $\vartheta$ as
%\begin{equation}
%  \label{eq:temperature-rescaling}
$
\frac{\vartheta}{\vartheta_{\mathrm{ref}}}
  =
  \left( \frac{\temp}{\tempref} \right)^{1 - m},
$
%\end{equation}
where $m \in (0, 1)$. By rescaling the temperature field one can identify the formula for the corresponding specific Helmholtz free energy---which will be different from the one given by \eqref{eq:fenergy}---and consequently, repeating the steps from Section \ref{sec:construction-of-the-functional}, one can obtain a whole family of functionals parameterized by $m$,
\begin{equation}
  \label{eq:functional-neq-m}
  \mathcal{V}_{\mathrm{neq}}^{m}
  \left(
    \left.
      \widetilde{\vec{W}}
    \right\|
    \widehat{\vec{W}}
  \right)
  =
  \int_{\Omega}
    \rho \cheatvolref \widehat{\temp}
    \left[
      \frac{\widetilde{\temp}}{\widehat{\temp}}
      -
      \frac{1}{m} 
      \left( 
        \left( 1 + \frac{\widetilde{\temp}}{\widehat{\temp}} \right)^m
        -
        1 
      \right)
    \right]
  \cvolumee
  +
  \int_{\Omega}
    \rho \fenergy_1 (\identity + \widetilde{\lcgGSII})
  \cvolumee
  +
  \int_{\Omega}
    \frac{1}{2} \rho \absnorm{\widetilde{\vecv}}^2
  \cvolumee.
\end{equation}
For any fixed $m \in (0, 1)$, functional $\mathcal{V}_{\mathrm{neq}}^{m}$ remains nonnegative and vanishes if and only if the perturbation $\widetilde{\vec{W}}$ vanishes.

For further reference, let us introduce the notation $\mathcal{V}_{\mathrm{th}}^{m}$ for the family of functionals
\begin{equation}
  \label{eq:functional-th-m}
  \mathcal{V}_{\mathrm{th}}^{m}
  \left(
    \left.
      \widetilde{\vec{W}}
    \right\|
    \widehat{\vec{W}}
  \right)
  =_{\bydefinition}
  \int_{\Omega}
    \rho \cheatvolref \widehat{\temp}
    \left[
      \frac{\widetilde{\temp}}{\widehat{\temp}}
      -
      \frac{1}{m} 
      \left( 
        \left( 1 + \frac{\widetilde{\temp}}{\widehat{\temp}} \right)^m
        -
        1 
      \right)
    \right]
  \cvolumee,
\end{equation}
that correspond to the thermal parts of functionals $\mathcal{V}_{\mathrm{neq}}^{m}$. \cite{dostalk.m.prusa.v.rajagopal.kr:unconditional} have shown that the time derivative of $\mathcal{V}_{\mathrm{th}}^{m}$ reads
\begin{multline}
  \label{eq:time-derivative-functional-th-m}
  \dd{\mathcal{V}_{\mathrm{th}}^{m}}{t}
  \left(
    \left.
      \widetilde{\vec{W}}
    \right\|
    \widehat{\vec{W}}
  \right)
  =
  -
  \int_{\Omega}
    4
    \frac{1 - m}{m^2}
    \kapparef \widehat{\temp} 
    \vectordot{
      \nabla 
      \left[ 
        \left( 
          1 + \frac{\widetilde{\temp}}{\widehat{\temp}} 
        \right)^{\frac{m}{2}} - 1 
      \right]
    }{
      \nabla 
      \left[ 
        \left( 
          1 + \frac{\widetilde{\temp}}{\widehat{\temp}} 
        \right)^{\frac{m}{2}} - 1 
      \right]
    }
  \, \cvolumee
  \\
  -
  \int_{\Omega}
    \frac{1 - m}{m}
    \rho \cheatvolref
    \left[
      \left( 1 + \frac{\widetilde{\temp}}{\widehat{\temp}} \right)^m - 1
    \right]
    \left(
      \vectordot{\widetilde{\vecv}}{\nabla \widehat{\temp}}  
    \right)
  \, \cvolumee
  \\
  +
  \int_{\Omega}
    \left(
      1 - \frac{1}{\left( 1 + \frac{\widetilde{\temp}}{\widehat{\temp}} \right)^{1 - m}}
    \right)
    \entprodctemp_{\mathrm{mech}}(\widehat{\vec{W}} + \widetilde{\vec{W}})
  \, \cvolumee,
\end{multline}
where $\entprodctemp_{\mathrm{mech}}(\widehat{\vec{W}} + \widetilde{\vec{W}})$ denotes the mechanical part of the entropy production. (\cite{dostalk.m.prusa.v.rajagopal.kr:unconditional} have shown~\eqref{eq:time-derivative-functional-th-m} for the Navier--Stokes--Fourier fluid with the entropy production term $\entprodctemp_{\mathrm{mech}} = 2 \nu \tensordot{\gradsym}{\gradsym}$. However, all the algebraic manipulations in~\cite{dostalk.m.prusa.v.rajagopal.kr:unconditional} hold also for more general entropy production term $\entprodctemp_{\mathrm{mech}}$.)

\section{Stability of the non-equilibrium steady state}
\label{sec:stability-of-the-non-equilibrium-steady-state}

A brief inspection of the right-hand side of~\eqref{eq:time-derivative-functional} reveals that the term with \emph{a~priori} unknown sign, that is the term
\begin{equation}
  \label{eq:4}
    \int_{\Omega}
    \rho \cheatvolref \ln \left( 1 + \frac{\widetilde{\temp}}{\widehat{\temp}} \right)
    \left(
      \vectordot{\widetilde{\vecv}}{\nabla \widehat{\temp}}  
    \right) 
    \, \cvolumee
    ,
\end{equation}
might be shown to be negligible provided that the velocity perturbation $\widetilde{\vec{v}}$ decays in time. This property is easy to show for the Navier--Stokes--Fourier fluid, see~\cite{dostalk.m.prusa.v.rajagopal.kr:unconditional} for details. \emph{Our objective is to recover the same property for the considered class of viscoelastic models.} This piece of information can be obtained by the analysis of the mechanical part $\mathcal{V}_{\mathrm{mech}}$ of the functional $\mathcal{V}_{\mathrm{neq}}$, see Section~\ref{sec:mechanical-quantities}.

Once we show that the norm of velocity perturbation is bounded by an exponentially decaying function, we can focus on the temperature perturbation only. Regarding the temperature perturbation, it is however straightforward to reuse results by~\cite{dostalk.m.prusa.v.rajagopal.kr:unconditional}. This is done in Section~\ref{sec:temperature-field}.

\subsection{Decay of perturbations -- mechanical quantities}
\label{sec:mechanical-quantities}
The formula~\eqref{eq:time-derivative-mechanical-functional} for the time derivative of the functional $\mathcal{V}_{\mathrm{mech}}$ can be rewritten explicitly as
\begin{multline}
  \label{eq:time-derivative-mechanical-functional-explicit}
  \dd{}{t}
  \int_{\Omega}
    \left(
      \rho \fenergy_1 (\identity + \widetilde{\lcgGSII})
      +
      \frac{1}{2} \rho \absnorm{\widetilde{\vecv}}^2
    \right)
  \cvolumee
  =
  -
  \int_{\Omega}
    2 \nu(\widehat{\temp} + \widetilde{\temp}) \tensordot{\widetilde{\gradsym}}{\widetilde{\gradsym}}
  \, \cvolumee
  \\
  -
  \int_{\Omega}
    \rho \frac{\mu}{\nu_1 (\widehat{\temp} + \widetilde{\temp})}
    \tensordot{
      \pd{\fenergy_1}{\lcgGSII} (\identity + \widetilde{\lcgGSII})
    }{
      \tensorf{f} (\identity + \widetilde{\lcgGSII})
    }
  \, \cvolumee.
\end{multline}
Since $\widetilde{\vecv}$ vanishes on the boundary, the Korn equality and the Poincar\'e inequality imply
\begin{equation}
  \label{eq:korn-poincare-inequality}
  \frac{1}{C_P} \,
  {\norm{\widetilde{\vecv}}}_{\sleb{2}{\Omega}}^2
  \le
  \int_{\Omega}
    2 \, \tensordot{\widetilde{\gradsym}}{\widetilde{\gradsym}}
  \, \cvolumee.
\end{equation}
Moreover, assumption \eqref{eq:fenergy-1-f-stability-inequality} gives us
\begin{equation}
  \label{eq:fenergy-1-f-stability-inequality-again}
  \fenergy_1 (\identity + \widetilde{\lcgGSII})
  \le
  C_{\tensorf{f}}
  \tensordot{\pd{\fenergy_1}{\lcgGSII} (\identity + \widetilde{\lcgGSII})}{\tensorf{f} (\identity + \widetilde{\lcgGSII})},
\end{equation}
where $C_{\tensorf{f}}$ is a positive constant dependent on the choice of $\tensorf{f}$.

Using inequalities \eqref{eq:korn-poincare-inequality}, \eqref{eq:fenergy-1-f-stability-inequality-again}, and boundedness of $\nu$ and $\nu_1$ from below and above respectively, we thus arrive at
\begin{multline}
  \label{eq:time-derivative-mechanical-functional-explicit-inequality}
  \dd{}{t}
  \int_{\Omega}
    \left(
      \rho \fenergy_1 (\identity + \widetilde{\lcgGSII})
      +
      \frac{1}{2} \rho \absnorm{\widetilde{\vecv}}^2
    \right)
  \cvolumee
  \le
  -
  \frac{2 \min_{s \in \R^+} \nu(s)}{\rho C_P} 
  \int_{\Omega}
    \frac{1}{2} \rho \absnorm{\widetilde{\vecv}}^2
  \cvolumee
  \\
  -
  \frac{\mu}{C_{\tensorf{f}} \max_{s \in \R^+} \nu_1(s)}
  \int_{\Omega}  
    \rho \fenergy_1 (\identity + \widetilde{\lcgGSII})
  \cvolumee.
\end{multline}
Consequently, estimate \eqref{eq:time-derivative-mechanical-functional-explicit-inequality} yields the following inequality for the time derivative of the functional $\mathcal{V}_{\mathrm{mech}}$
\begin{equation}
  \label{eq:time-derivative-mechanical-functional-final-inequality}
  \dd{\mathcal{V}_{\mathrm{mech}}}{t}
  \left(
    \left.
      \widetilde{\vec{W}}
    \right\|
    \widehat{\vec{W}}
  \right)
  \le
  -
  \, C_{\mathrm{mech}}
  \, \mathcal{V}_{\mathrm{mech}}
  \left(
    \left.
      \widetilde{\vec{W}}
    \right\|
    \widehat{\vec{W}}
  \right),
\end{equation}
where we have denoted
\begin{equation}
  \label{eq:c-mechanical}
  C_{\mathrm{mech}}
  =_{\bydefinition}
  \min
  \left\{
    \frac{2 \min_{s \in \R^+} \nu(s)}{\rho C_P},
    \,
    \frac{\mu}{C_{\tensorf{f}} \max_{s \in \R^+} \nu_1(s)}
  \right\}.
\end{equation}
It then follows that
\begin{equation}
  \label{eq:mechanical-functional-exp-decay}
  \mathcal{V}_{\mathrm{mech}}
  \left(
    \left.
      \widetilde{\vec{W}}
    \right\|
    \widehat{\vec{W}}
  \right)
  \le
  \left.
    \mathcal{V}_{\mathrm{mech}}
    \left(
      \left.
        \widetilde{\vec{W}}
      \right\|
      \widehat{\vec{W}}
    \right)
  \right|_{t=0}
  e^{-C_{\mathrm{mech}} t},
\end{equation}
which further implies
\begin{subequations}
  \label{eq:exp-decays}
  \begin{align}
    \label{eq:exp-decays-velocity-norm}
    {\norm{\widetilde{\vecv}}}_{\sleb{2}{\Omega}}^2
    &\le
    \frac{2}{\rho}
    \left.
      \mathcal{V}_{\mathrm{mech}}
      \left(
        \left.
          \widetilde{\vec{W}}
        \right\|
        \widehat{\vec{W}}
      \right)
    \right|_{t=0}
    e^{-C_{\mathrm{mech}} t},
    \\
    \label{eq:exp-decays-fenergy-1}
    \int_{\Omega}
      \fenergy_1 (\identity + \widetilde{\lcgGSII})
    \, \cvolumee
    &\le
    \frac{1}{\rho}
    \left.
      \mathcal{V}_{\mathrm{mech}}
      \left(
        \left.
          \widetilde{\vec{W}}
        \right\|
        \widehat{\vec{W}}
      \right)
    \right|_{t=0}
    e^{-C_{\mathrm{mech}} t}.
  \end{align}
\end{subequations}
Estimates \eqref{eq:exp-decays} yield the desired result. The perturbations $\widetilde{\vecv}$ and $\widetilde{\lcgGSII}$ vanish as time goes to infinity. (See the assumption \eqref{eq:fenergy-1-assumption-zero-value}.) Note also that~\eqref{eq:exp-decays-fenergy-1} implies only the decay of quantity
$    \int_{\Omega}
      \fenergy_1 (\identity + \widetilde{\lcgGSII})
    \, \cvolumee
$, while this quantity might be difficult to interpret as a convergence of $\widetilde{\lcgGSII}$ to zero in a norm. Still there is a relation between this quantity and a reasonable metric on the set of spatially distributed symmetric positive definite matrices. (The metric is constructed using the Bures--Wasserstein distance on the set of positive definite matrices, see~\cite{bhatia.r.jain.t.ea:on}.) For details regarding this concept we refer the interested reader to~\cite{dostal.m.prusa.v.ea:finite}.

\subsection{Decay of perturbation -- temperature}
\label{sec:temperature-field}
Having obtained an upper bound on the norm of the velocity perturbation, we reuse the results by~\cite{dostalk.m.prusa.v.rajagopal.kr:unconditional} for the standard Navier--Stokes--Fourier fluid occupying a mechanically isolated vessel with spatially non-uniform wall temperature. The authors show that the spatially inhomogeneous steady temperature field $\widehat{\temp}$ is stable irrespective of the initial temperature field. The derivation rests upon the usage of the family of functionals $\mathcal{V}_{\mathrm{th}}^m$ introduced in Section~\ref{sec:family-of-functionals-th} and exploits the fact that ${\norm{\widetilde{\vecv}}}_{\sleb{2}{\Omega}}$ is bounded from above by an exponentially decaying function. Further, the entropy production $\entprodctemp_{\mathrm{mech}}(\widehat{\vec{W}} + \widetilde{\vec{W}})$ must be a nonnegative quantity that vanishes at equilibrium. Since these properties hold in our case as well, see \eqref{eq:exp-decays-velocity-norm}, we can directly generalise the result of \cite{dostalk.m.prusa.v.rajagopal.kr:unconditional} to the viscoelastic models (Oldroyd-B, Giesekus, FENE-P, Johnson--Segalman, Phan-Thien--Tanner) described in~\ref{sec:admissible-viscoelastic-models}.

In particular, one can show that for $n, m \in (0, 1)$, $n > m > \frac{n}{2}$ the functional 
\begin{equation}
  \label{eq:y-m-n}
  \mathcal{Y}_{\mathrm{th}}^{m,n}
  \left(
    \left.
      \widetilde{\vec{W}}
    \right\|
    \widehat{\vec{W}}
  \right)
  =_{\bydefinition}
  \mathcal{V}_{\mathrm{th}}^{m}
  \left(
    \left.
      \widetilde{\vec{W}}
    \right\|
    \widehat{\vec{W}}
  \right)
  -
  \mathcal{V}_{\mathrm{th}}^{n}
  \left(
    \left.
      \widetilde{\vec{W}}
    \right\|
    \widehat{\vec{W}}
  \right)
\end{equation}
decays to zero as time goes to infinity. Specifically, according to the definition \eqref{eq:functional-th-m}, this translates to 
\begin{equation}
  \label{eq:temperature-vanishes}
  \int_{\Omega}
    \rho \cheatvolref \widehat{\temp}
    \left[
      \frac{1}{n} \left( 1 + \frac{\widetilde{\temp}}{\widehat{\temp}} \right)^n
      -
      \frac{1}{m} \left( 1 + \frac{\widetilde{\temp}}{\widehat{\temp}} \right)^m
      +
      \frac{n - m}{mn}
    \right]
  \cvolumee
  \xrightarrow{t \rightarrow +\infty}
  0.
\end{equation}
Using \cite[Corollary 1]{dostalk.m.prusa.v.rajagopal.kr:unconditional} we also see that~\eqref{eq:temperature-vanishes} implies the decay of the relative entropy in any Lebesgue space $\sleb{p}{\Omega}$, $p \in [1, +\infty)$.

In order to obtain \eqref{eq:temperature-vanishes} one needs to show that all the terms on the right-hand side of \eqref{eq:time-derivative-functional-th-m} are finite if we integrate them with respect to time from zero to infinity. This is where \eqref{eq:exp-decays-velocity-norm} comes into play. Finally, the convergence result \eqref{eq:temperature-vanishes} then follows from a lemma on the decay of integrable functions, see \cite[Lemma 1.2]{zheng.s:asymptotic}, applied to the functional $\mathcal{Y}_{\mathrm{th}}^{m,n}$, see~\cite{dostalk.m.prusa.v.rajagopal.kr:unconditional} for details.  

\section{Conclusion}
\label{sec:conclusion}
We have investigated the stability of a spatially inhomogeneous non-equilibrium steady state in a thermodynamically open system. Specifically, we have dealt with an incompressible heat conducting viscoelastic fluid occupying a vessel with spatially non-uniform wall temperature. The steady state in this system is characterised by the zero velocity field $\vec{v}$ and a trivial~$\lcgGSII$ field, while the temperature field $\temp$ is the solution of the steady heat equation.

Assuming that the governing equations possess the~\emph{classical solution} that exists for all times, we have shown that the steady state is stable irrespective of the initial conditions and of the shape of the vessel. (The perturbations decay to zero as time goes to infinity.) We have thus generalised the results by~\cite{dostalk.m.prusa.v.rajagopal.kr:unconditional}, who have investigated the same stability problem for the incompressible Navier--Stokes--Fourier fluid. Our analysis is general enough to capture a wide range of viscoelastic models including the Oldroyd-B model, the Giesekus model, the FENE-P model, the Johnson--Segalman model, and the Phan--Thien--Tanner model. 

\bibliography{vit-prusa,heat-stability-viscoelastic}

%\newpage
\appendix

\section{Admissible viscoelastic models}
\label{sec:admissible-viscoelastic-models}
Let us show that the Oldroyd-B model, the Giesekus model, the FENE-P model, the Johnson--Segalman model, and the Phan--Thien--Tanner model satisfy the structural assumptions introduced in Section~\ref{sec:helmholtz-free-energy} and Section~\ref{sec:entropy-production}. In particular, we show that all the assumptions imposed on the scalar function $\fenergy_1$ (which determines the specific Helmholtz energy) and the tensorial function $\tensorf{f}$ (which determines the specific entropy production) are fulfilled. Let us reiterate the requirements from Section~\ref{sec:helmholtz-free-energy} and Section~\ref{sec:entropy-production} here.

\emph{First}, for a given $\fenergy_1 : \R_{>}^{3 \times 3} \to \R$, where $\R_{>}^{3 \times 3}$ denotes the set of symmetric positive definite~$3 \times 3$ matrices, we need to verify that
\begin{subequations}
  \label{eq:fenergy-1-assumptions-appendix}
  \begin{gather}
    \label{eq:fenergy-1-assumption-zero-value-appendix}
    \fenergy_1 (\lcgGSII) \ge 0, \qquad \fenergy_1(\lcgGSII) = 0 \iff \lcgGSII = \identity,
    \\
    \label{eq:fenergy-1-assumption-derivative-zero-value-appendix}
    \pd{\fenergy_1}{\lcgGSII}(\lcgGSII) = \tensorq{0} \iff \lcgGSII = \identity,
    \\
    \label{eq:fenergy-1-assumption-commutativity-appendix}
    \lcgGSII \pd{\fenergy_1}{\lcgGSII}(\lcgGSII) = \pd{\fenergy_1}{\lcgGSII}(\lcgGSII) \lcgGSII,
  \end{gather}
\end{subequations}
hold for any symmetric positive definite tensor $\lcgGSII$.

\emph{Second}, a given tensorial function $\tensorf{f} : \R_{>}^{3 \times 3} \to \R_{>}^{3 \times 3}$ must, for any symmetric positive definite tensor $\lcgGSII$, meet the following requirements
\begin{subequations}
  \label{eq:f-assumptions-appendix}
  \begin{gather}
    \label{eq:f-zero-appendix}
    \tensorf{f}(\lcgGSII) = \tensorq{0}
    \iff
    \lcgGSII = \identity,
    \\
    \label{eq:pd-fenergy-1-f-nonnegative-appendix}
    \tensordot{\pd{\fenergy_1}{\lcgGSII} (\lcgGSII)}{\tensorf{f}(\lcgGSII)}
    \ge
    0,
    \\
    \label{eq:fenergy-1-f-stability-inequality-appendix}
    \fenergy_1 (\lcgGSII)
    \le
    C_{\tensorf{f}}
    \tensordot{\pd{\fenergy_1}{\lcgGSII} (\lcgGSII)}{\tensorf{f} (\lcgGSII)},
  \end{gather}
\end{subequations}
where $C_{\tensorf{f}}$ is a positive constant dependent on the choice of $\tensorf{f}$. The last assumption \eqref{eq:fenergy-1-f-stability-inequality-appendix} is crucial for obtaining the stability result~\eqref{eq:exp-decays}.

\subsection{Oldroyd-B model}
\label{sec:maxwell-oldroyd-b-model}
As discussed in Section \ref{sec:kinematics} we set $a = 1$ for the Oldroyd-B model and use the notation $\lcgII$ instead of $\lcgGSII$ for the additional tensorial quantity in the Cauchy stress tensor.

%\subsubsection{Helmholtz free energy $\fenergy_1$}
%\label{sec:maxwell-oldroyd-b-helmholtz-free-energy}
The ``elastic'' part $\fenergy_1$ of the specific free energy for the Oldroyd-B model reads
\begin{equation}
  \label{eq:fenergy-1-oldroyd-b}
  \fenergy_1 (\lcgII)
  =_{\bydefinition}
  \frac{\mu}{2 \rho}
  \left(
    \Tr \lcgII - 3 - \ln \det \lcgII
  \right).
\end{equation}
Using the identity $\ln \det \lcgII = \Tr \ln \lcgII$ we can write
\begin{equation}
  \label{eq:fenergy-1-oldroyd-b-nonnegativity}
  \fenergy_1 (\lcgII)
  =
  \frac{\mu}{2 \rho} \Tr \left( \lcgII - \identity - \ln \lcgII \right)
  =
  \frac{\mu}{2 \rho}
  \sum_{i=1}^{3} \left( \lambda_i - 1 - \ln \lambda_i \right) ,
\end{equation}
where $\{\lambda_i\}_{i=1}^3$ denote eigenvalues of the symmetric positive  definite tensor $\lcgII$. Since the function $f(x) =_{\bydefinition} x - 1 - \ln x$ is nonnegative for $x > 0$ and vanishes if and only if $x = 1$, we obtain the validity of \eqref{eq:fenergy-1-assumption-zero-value-appendix}.

The derivative of $\fenergy_1$ with respect to $\lcgII$ reads
\begin{equation}
  \label{eq:fenergy-1-derivative-oldroyd-b}
  \pd{\fenergy_1}{\lcgII} (\lcgII)
  =
  \frac{\mu}{2 \rho}
  \left(
    \identity - \inverse{\lcgII}
  \right),
\end{equation}
and we immediately see that assumptions \eqref{eq:fenergy-1-assumption-derivative-zero-value-appendix} and \eqref{eq:fenergy-1-assumption-commutativity-appendix} are both fulfilled.

%\subsubsection{Tensorial function $\tensorf{f}$}
%\label{sec:maxwell-oldroyd-b-tensorial-function}
The tensorial function $\tensorf{f}$ for the Oldroyd-B model reads
\begin{equation}
  \label{eq:f-oldroyd-b}
  \tensorf{f}(\lcgII) =_{\bydefinition} \lcgII - \identity.
\end{equation}
The requirement \eqref{eq:f-zero-appendix} is obviously satisfied. To verify the validity of \eqref{eq:pd-fenergy-1-f-nonnegative-appendix} let us write
\begin{equation}
  \label{eq:f-assumption-1-oldroyd-b}
  \tensordot{\pd{\fenergy_1}{\lcgII}}{\tensorf{f}(\lcgII)}
  =
  \frac{\mu}{2 \rho}
  \Tr
  \left(
    \lcgII - 2 \identity + \inverse{\lcgII}
  \right)
  =
  \frac{\mu}{2 \rho}
  \sum_{i=1}^3
  \left(
    \lambda_i - 2 + \frac{1}{\lambda_i}
  \right).
\end{equation}
Since the function $g(x) =_{\bydefinition} x - 2 + 1/x$ is nonnegative for $x > 0$ and vanishes if and only if $x = 1$, we see that \eqref{eq:pd-fenergy-1-f-nonnegative-appendix} is fulfilled. 

Finally, we want to show that the inequality \eqref{eq:fenergy-1-f-stability-inequality-appendix} holds, which for the given $\fenergy_1$ and $\tensorf{f}$ translates to
\begin{equation}
  \label{eq:1-oldroyd-b}
  \Tr \lcgII - 3 - \ln \det \lcgII
  \le
  C_{\tensorf{f}} \Tr \left( \lcgII - 2 \identity + \inverse{\lcgII} \right).
\end{equation}
Taking $C_{\tensorf{f}} =_{\bydefinition} 1$ and using the identity $\ln \det \lcgII = \Tr \ln \lcgII$ we can rewrite \eqref{eq:1-oldroyd-b} as
\begin{equation}
  \label{eq:2-oldroyd-b}
  0 
  \le 
  \Tr \left( \inverse{\lcgII} + \ln \lcgII - \identity \right).
\end{equation}
But the right-hand side of \eqref{eq:2-oldroyd-b} is indeed nonnegative since
\begin{equation}
  \label{eq:3-oldroyd-b}
  \Tr \left( \inverse{\lcgII} + \ln \lcgII - \identity \right)
  =
  \sum_{i=1}^3 \left( \frac{1}{\lambda_i} + \ln \lambda_i - 1 \right),
\end{equation}
and the function $h(x) =_{\bydefinition} 1/x + \ln x - 1$ is nonnegative for $x > 0$ and vanishes if and only if $x = 1$. One can easily show that $C_{\tensorf{f}} = 1$ is optimal, that is taking $C_{\tensorf{f}}$ smaller would violate \eqref{eq:fenergy-1-f-stability-inequality-appendix}.

\subsection{Giesekus model}
\label{sec:giesekus-model}
As discussed in Section \ref{sec:kinematics} we set $a = 1$ for the Giesekus model and use the notation $\lcgII$ instead of $\lcgGSII$ for the additional tensorial quantity in the Cauchy stress tensor.

%\subsubsection{Helmholtz free energy $\fenergy_1$}
%\label{sec:giesekus-helmholtz-free-energy}
The ``elastic'' part $\fenergy_1$ of the specific free energy for the Giesekus model reads
\begin{equation}
  \label{eq:fenergy-1-giesekus}
  \fenergy_1 (\lcgII)
  =_{\bydefinition}
  \frac{\mu}{2 \rho}
  \left(
    \Tr \lcgII - 3 - \ln \det \lcgII
  \right).
\end{equation}
This is the same specific free energy as in the case of Oldroyd-B model, and we already know that this choice of Helmholtz free energy satisfies~\eqref{eq:fenergy-1-assumptions-appendix}, see~\ref{sec:maxwell-oldroyd-b-model}.

%\subsubsection{Tensorial function $\tensorf{f}$}
%\label{sec:giesekus-tensorial-function}
The tensorial function $\tensorf{f}$ for the Giesekus model reads
\begin{equation}
  \label{eq:f-giesekus}
  \tensorf{f}(\lcgII) 
  =_{\bydefinition} 
  \alpha \lcgII^2 + (1 - 2 \alpha) \lcgII - (1 - \alpha) \identity,
\end{equation}
where $\alpha \in (0, 1)$ is a model parameter. Since $\lcgII$ is a symmetric positive definite tensor it is diagonalizable and we can thus easily show that \eqref{eq:f-zero-appendix} indeed holds. To verify assumption \eqref{eq:pd-fenergy-1-f-nonnegative-appendix} let us write
\begin{multline}
  \label{eq:f-assumption-1-giesekus}
  \tensordot{\pd{\fenergy_1}{\lcgII}}{\tensorf{f}(\lcgII)}
  =
  \frac{\mu}{2 \rho}
  \Tr
  \left[
    \alpha \lcgII^2 + (1 - 3 \alpha) \lcgII - (2 - 3 \alpha) \identity
    + (1 - \alpha) \inverse{\lcgII}
  \right]
  \\
  =
  \frac{\mu}{2 \rho}
  \sum_{i=1}^3
  \left(
    \alpha \lambda_i^2 + (1 - 3 \alpha) \lambda_i - (2 - 3 \alpha) + (1 - \alpha) \frac{1}{\lambda_i}
  \right).
\end{multline}
It is straightforward to show that for $\alpha \in (0, 1)$ the function $g_{\alpha}(x) =_{\bydefinition} \alpha x^2 + (1 - 3 \alpha) x - (2 - 3 \alpha) + (1 - \alpha) 1/x$ is nonnegative for $x > 0$, and that it vanishes if and only if $x = 1$. Assumption \eqref{eq:pd-fenergy-1-f-nonnegative-appendix} is thus fulfilled. Finally, we want to show that the inequality \eqref{eq:fenergy-1-f-stability-inequality-appendix} holds, which for the given $\fenergy_1$ and $\tensorf{f}$ translates to
\begin{equation}
  \label{eq:1-giesekus}
  \Tr \lcgII - 3 - \ln \det \lcgII
  \le
  C_{\tensorf{f}} 
  \Tr
  \left[
    \alpha \lcgII^2 + (1 - 3 \alpha) \lcgII - (2 - 3 \alpha) \identity
    + (1 - \alpha) \inverse{\lcgII}
  \right].
\end{equation}
Taking $C_{\tensorf{f}} =_{\bydefinition} \frac{1}{1 - \alpha}$ and using the identity $\ln \det \lcgII = \Tr \ln \lcgII$, a simple manipulation reveals that \eqref{eq:1-giesekus} is equivalent to
\begin{equation}
  \label{eq:2-giesekus}
  0 
  \le 
  \frac{\alpha}{1 - \alpha}
  \Tr \left[ \left( \lcgII - \identity \right)^2 \right]
  +
  \Tr \left( \inverse{\lcgII} + \ln \lcgII - \identity \right).
\end{equation}
However, the first term on the right-hand side of~\eqref{eq:2-giesekus} is obviously nonnegative and the second term is nonnegative as well as has been shown in the case of the Oldroyd-B model, see~\eqref{eq:3-oldroyd-b}.  

\subsection{FENE-P model}
\label{sec:fene-p-model}
As discussed in Section \ref{sec:kinematics} we set $a = 1$ for the FENE-P model and use the notation $\lcgII$ instead of $\lcgGSII$ for the additional tensorial quantity in the Cauchy stress tensor.

%\subsubsection{Helmholtz free energy $\fenergy_1$}
%\label{sec:fene-p-helmholtz-free-energy}
The ``elastic'' part $\fenergy_1$ of the specific free energy for the FENE-P model reads 
\begin{equation}
  \label{eq:fenergy-1-fene-p}
  \fenergy_1 (\lcgII)
  =_{\bydefinition}
  \frac{\mu}{2 \rho}
  \left[
    -
    b \ln \left( 1 - \frac{1}{b} \Tr \lcgII \right)
    +
    b \ln \left( 1 - \frac{3}{b} \right)
    -
    \inverse{\left( 1 - \frac{3}{b} \right)} \ln \det \lcgII
  \right],
\end{equation}
where $b > 3$ is a model parameter and $\Tr \lcgII < b$. We want to show that
\begin{equation}
  \label{eq:0a-fene-p}
  -
  b \ln \left( 1 - \frac{1}{b} \Tr \lcgII \right)
  +
  b \ln \left( 1 - \frac{3}{b} \right)
  -
  \inverse{\left( 1 - \frac{3}{b} \right)} \ln \det \lcgII
  \ge
  0.
\end{equation}
Inequality \eqref{eq:0a-fene-p} can be rewritten in the following form
\begin{equation}
  \label{eq:0b-fene-p}
  b \ln 
  \left[ 
    \frac{b - 3}{b - \Tr \lcgII} \left( \det \lcgII \right)^{\frac{1}{3-b}}
  \right]
  \ge
  0,
\end{equation}
and thus, it suffices to investigate whether
\begin{equation}
  \label{eq:0c-fene-p}
  \frac{b - 3}{b - \Tr \lcgII} \left( \det \lcgII \right)^{\frac{1}{3-b}}
  \ge
  1,
\end{equation}
holds. Since $\lcgII$ is symmetric positive definite, the standard inequality of arithmetic and geometric means yields $\Tr \lcgII \ge 3 (\det \lcgII)^{\frac{1}{3}}$. Consequently, it suffices to prove the following inequality
\begin{equation}
  \label{eq:0d-fene-p}
  \frac{b - 3}{b - 3 (\det \lcgII)^{\frac{1}{3}}} \left( \det \lcgII \right)^{\frac{1}{3-b}}
  \ge
  1,
\end{equation}
which can be further rewritten as
\begin{equation}
  \label{eq:0e-fene-p}
  3 \left[ \left( \det \lcgII \right)^{\frac{1}{3}} - 1 \right]
  -
  (3 - b) \left[ \left( \det \lcgII \right)^{\frac{1}{3 - b}} - 1 \right]
  \ge
  0.
\end{equation}
A simple analysis reveals that the function $f_{(r, s)}(x) =_{\bydefinition} r (x^{\frac{1}{r}} - 1) - s (x^{\frac{1}{s}} - 1)$, where $r > 0$, $s < 0$, is nonnegative for $x > 0$ and vanishes if and only if $x = 1$. We have thus proved that $\psi (\lcgII) \ge 0$, and that $\psi (\lcgII) = 0$ implies $\det \lcgII = 1$. It is then straightforward to check that $\psi (\lcgII) = 0$, if and only if $\lcgII = \identity$, and the verification of assumption \eqref{eq:fenergy-1-assumption-zero-value-appendix} is thus complete. % Indeed, once we know that $\psi (\lcgII) = 0$, then we know that $\det \lcgII = 1$ and that $\Tr \lcgII = 3$. The inequality of arithmetic and geometric means implies that $\Tr \lcgII \geq 3 \left(\det \lcgII\right)^{\frac{1}{3}}$, and we know that we have equality here. However, the AM-AG inequality is an equality if and only if all the elements are equal to zero. 

The derivative of $\fenergy_1$ with respect to $\lcgII$ reads
\begin{equation}
  \label{eq:fenergy-1-derivative-fene-p}
  \pd{\fenergy_1}{\lcgII} (\lcgII)
  =
  \frac{\mu}{2 \rho}
  \left[
    \inverse{\left( 1 - \frac{1}{b} \Tr \lcgII \right)}\identity
    -
    \inverse{\left( 1 - \frac{3}{b} \right)} \inverse{\lcgII}
  \right],
\end{equation}
and we immediately see that assumption \eqref{eq:fenergy-1-assumption-commutativity-appendix} is fulfilled. Further, the fact that $\lcgII$ is diagonalizable yields the validity of \eqref{eq:fenergy-1-assumption-derivative-zero-value-appendix}.

%\subsubsection{Tensorial function $\tensorf{f}$}
%\label{sec:fene-p-tensorial-function}
The tensorial function $\tensorf{f}$ for the FENE-P model reads
\begin{equation}
  \label{eq:f-fene-p}
  \tensorf{f} (\lcgII)
  =_{\bydefinition}
  \inverse{\left( 1 - \frac{1}{b} \Tr \lcgII \right)} \lcgII
  -
  \inverse{\left( 1 - \frac{3}{b} \right)} \identity.
\end{equation}
Diagonalization of the tensor $\lcgII$ can be used to confirm the validity of requirement \eqref{eq:f-zero-appendix}. To verify assumption \eqref{eq:pd-fenergy-1-f-nonnegative-appendix} let us write
\begin{equation}
  \label{eq:1-fene-p}
  \tensordot{\pd{\fenergy_1}{\lcgII}}{\tensorf{f}(\lcgII)}
  =
  \frac{\mu}{2 \rho}
  \Tr
  \left[
    \left( 1 - \frac{1}{b} \Tr \lcgII \right)^{-2} \lcgII
    -
    2 \inverse{\left( 1 - \frac{1}{b} \Tr \lcgII \right)} \inverse{\left( 1 - \frac{3}{b} \right)} \identity
    +
    \left( 1 - \frac{3}{b} \right)^{-2} \inverse{\lcgII}
  \right].
\end{equation}
The right-hand side of \eqref{eq:1-fene-p} can be rewritten using the eigenvalues of $\lcgII$ as
\begin{multline}
  \label{eq:2-fene-p}
  \frac{\mu}{2 \rho}
  \sum_{i=1}^3
  \left[
    \left( 1 - \frac{1}{b} \Tr \lcgII \right)^{-2} \lambda_i 
    - 
    2 \inverse{\left( 1 - \frac{1}{b} \Tr \lcgII \right)} \inverse{\left( 1 - \frac{3}{b} \right)}
    + 
    \left( 1 - \frac{3}{b} \right)^{-2} \frac{1}{\lambda_i}
  \right]
  \\
  =
  \frac{\mu}{2 \rho}
  \sum_{i=1}^3
  \frac{1}{\lambda_i}
  \left[
    \left( 1 - \frac{1}{b} \Tr \lcgII \right)^{-1} \lambda_i 
    -
    \left( 1 - \frac{3}{b} \right)^{-1}
  \right]^2,
\end{multline}
and we immediately see that the right-hand side of \eqref{eq:2-fene-p} is nonnegative.

Finally, we want to show that the inequality \eqref{eq:fenergy-1-f-stability-inequality-appendix} holds, which for the given $\fenergy_1$ and $\tensorf{f}$ translates to
\begin{multline}
  \label{eq:3-fene-p}
  -
  b \ln \left( 1 - \frac{1}{b} \Tr \lcgII \right)
  +
  b \ln \left( 1 - \frac{3}{b} \right)
  -
  \inverse{\left( 1 - \frac{3}{b} \right)} \ln \det \lcgII
  \\
  \le
  C_{\tensorf{f}}
  \Tr
  \left[
    \left( 1 - \frac{1}{b} \Tr \lcgII \right)^{-2} \lcgII
    -
    2 \inverse{\left( 1 - \frac{1}{b} \Tr \lcgII \right)} \inverse{\left( 1 - \frac{3}{b} \right)} \identity
    +
    \left( 1 - \frac{3}{b} \right)^{-2} \inverse{\lcgII}
  \right].
\end{multline}
Taking $C_{\tensorf{f}} =_{\bydefinition} 1 - \frac{3}{b}$ and using the identity $\ln \det \lcgII = \Tr \ln \lcgII$, a simple manipulation reveals that~\eqref{eq:3-fene-p} is equivalent to
\begin{multline}
  \label{eq:4-fene-p}
  \left( 1 - \frac{3}{b} \right) \left( 1 - \frac{1}{b} \Tr \lcgII \right)^{-2} \Tr \lcgII
  -
  6 \inverse{\left( 1 - \frac{1}{b} \Tr \lcgII \right)}
  +
  b \ln \left( 1 - \frac{1}{b} \Tr \lcgII \right)
  \\
  -
  b \ln \left( 1 - \frac{3}{b} \right)
  +
  3\inverse{\left( 1 - \frac{3}{b} \right)}
  +
  \inverse{\left( 1 - \frac{3}{b} \right)}
  \Tr 
  \left[
    \inverse{\lcgII} + \ln \lcgII - \identity 
  \right]
  \ge
  0.
\end{multline}
The last term on the left-hand side of \eqref{eq:4-fene-p} is nonnegative, see \eqref{eq:3-oldroyd-b}. Hence, it suffices to show that
\begin{multline}
  \label{eq:5-fene-p}
  \left( 1 - \frac{3}{b} \right) \left( 1 - \frac{1}{b} \Tr \lcgII \right)^{-2} \Tr \lcgII
  -
  6 \inverse{\left( 1 - \frac{1}{b} \Tr \lcgII \right)}
  \\
  +
  b \ln \left( 1 - \frac{1}{b} \Tr \lcgII \right)
  -
  b \ln \left( 1 - \frac{3}{b} \right)
  +
  3\inverse{\left( 1 - \frac{3}{b} \right)}
  \ge
  0.
\end{multline}
Since we know that $0 < \Tr \lcgII < b$, let us write $\Tr \lcgII = \varepsilon b$, where $\varepsilon \in (0, 1)$. Moreover, let us denote
\begin{equation}
  \label{eq:6-fene-p}
  f_b(\varepsilon)
  =_{\bydefinition}
  \left( 1 - \frac{3}{b} \right) \left( 1 - \varepsilon \right)^{-2} \varepsilon b
  -
  6 \inverse{\left( 1 - \varepsilon \right)}
  \\
  +
  b \ln \left( 1 - \varepsilon \right)
  -
  b \ln \left( 1 - \frac{3}{b} \right)
  +
  3\inverse{\left( 1 - \frac{3}{b} \right)}.
\end{equation}
Inequality \eqref{eq:5-fene-p} then transforms into the question whether the function $f_b(\varepsilon)$ is nonnegative for $\varepsilon \in (0,1)$ and $b > 3$. A tedious but straightforward analysis of $f_b$ reveals that this is indeed the case and assumption \eqref{eq:fenergy-1-f-stability-inequality-appendix} is thus verified.

\subsection{Johnson--Segalman model}
\label{sec:johnson-segalman-validity}

%\subsubsection{Helmholtz free energy $\fenergy_1$}
%\label{sec:js-helmholtz-free-energy}
The ``elastic'' part $\fenergy_1$ of the specific free energy for the Johnson--Segalman model reads
\begin{equation}
  \label{eq:fenergy-1-js}
  \fenergy_1 (\lcgGSII)
  =_{\bydefinition}
  \frac{\mu}{2 \rho}
  \left(
    \Tr \lcgGSII - 3 - \ln \det \lcgGSII
  \right).
\end{equation}
We see that apart from the usage of the tensorial quantity $\lcgGSII$ instead of $\lcgII$, the specific free energy of the Johnson--Segalman model is the same as of the Oldroyd-B model. Assumptions \eqref{eq:fenergy-1-assumptions-appendix} have thus been already verified in~\ref{sec:maxwell-oldroyd-b-model}.

%\subsubsection{Tensorial function $\tensorf{f}$}
%\label{sec:js-tensorial-function}
The tensorial function $\tensorf{f}$ for the Johnson--Segalman model reads
\begin{equation}
  \label{eq:f-js}
  \tensorf{f} (\lcgGSII)
  =_{\bydefinition}
  \lcgGSII - \identity,
\end{equation}
and again we see that assumptions~\eqref{eq:f-assumptions-appendix} have already been verified in~\ref{sec:maxwell-oldroyd-b-model} since the only difference between the tensorial function of the Johnson--Segalman model and of the Oldroyd-B model lies in the different physical interpretation of its tensorial argument.

\subsection{Phan--Thien--Tanner model}
\label{sec:phan-thien-tanner-model}

%\subsubsection{Helmholtz free energy $\fenergy_1$}
%\label{sec:ptt-helmholtz-free-energy}
The ``elastic'' part $\fenergy_1$ of the specific free energy for the Phan--Thien--Tanner model reads
\begin{equation}
  \label{eq:fenergy-1-ptt}
  \fenergy_1 (\lcgGSII)
  =_{\bydefinition}
  \frac{\mu}{2 \rho}
  \left(
    \Tr \lcgGSII - 3 - \ln \det \lcgGSII
  \right).
\end{equation}
The formula \eqref{eq:fenergy-1-ptt} is the same as for the Johnson-Segalman model and in turn as for the Oldroyd-B model apart from its different tensorial argument. Assumptions \eqref{eq:fenergy-1-assumptions-appendix} have thus been already verified in~\ref{sec:maxwell-oldroyd-b-model}.

%\subsubsection{Tensorial function $\tensorf{f}$}
%\label{sec:ptt-tensorial-function}
The tensorial function $\tensorf{f}$ for the Phan--Thien--Tanner model reads
\begin{equation}
  \label{eq:f-ptt}
  \tensorf{f} (\lcgGSII)
  =_{\bydefinition}
  \exponential{ p \Tr \left( \lcgGSII - \identity \right) }
  \left( \lcgGSII - \identity \right),
\end{equation}
where $p > 0$ is a model parameter.\footnote{
  For the sake of simplicity, we consider the \emph{exponential} Phan--Thien--Tanner model as proposed by \cite{phan-thien.n:non-linear}, which is the model given by the tensorial function~\eqref{eq:f-ptt}. There are other models referred to as the Phan--Thien--Tanner model, see for example the \emph{linear} Phan--Thien--Tanner model introduced in \cite{phan-thien.n.tanner.ri:new}. In this case the tensorial function $\tensorf{f}$ is given by by the formula
    \begin{equation*}
      \tensorf{f} (\lcgGSII)
      =_{\bydefinition}
      \left[ 1 + p \Tr \left( \lcgGSII - \identity \right) \right]
      \left( \lcgGSII - \identity \right),
    \end{equation*}
  where $p \in (0, 1/3]$. Assumptions \eqref{eq:f-assumptions-appendix} could be easily shown to hold for this case as well.
}
From \eqref{eq:f-ptt} it can be immediately seen that \eqref{eq:f-zero-appendix} is fulfilled. To verify the validity of \eqref{eq:pd-fenergy-1-f-nonnegative-appendix} let us write
\begin{equation}
  \label{eq:f-assumption-1-ptt}
  \tensordot{\pd{\fenergy_1}{\lcgGSII}}{\tensorf{f}(\lcgGSII)}
  =
  \frac{\mu}{2 \rho}
  \exponential{ p \Tr \left( \lcgGSII - \identity \right) }
  \Tr
  \left(
    \lcgGSII - 2 \identity + \inverse{\lcgGSII}
  \right).
\end{equation}
The nonnegativity of the right-hand side of \eqref{eq:f-assumption-1-ptt} can then be obtained just as in the case of the Oldroyd-B model, see \eqref{eq:f-assumption-1-oldroyd-b}.

It remains to verify assumption \eqref{eq:fenergy-1-f-stability-inequality-appendix} which for the given $\fenergy_1$ and $\tensorf{f}$ translates to
\begin{equation}
  \label{eq:1-ptt}
  \left(
    \Tr \lcgGSII - 3 - \ln \det \lcgGSII
  \right)
  \le
  C_{\tensorf{f}}
  \exponential{ p \Tr \left( \lcgGSII - \identity \right) }
  \Tr
  \left(
    \lcgGSII - 2 \identity + \inverse{\lcgGSII}
  \right).
\end{equation}
Taking $C_{\tensorf{f}} =_{\bydefinition} \exponential{3p}$ and using the identity $\ln \det \lcgGSII = \Tr \ln \lcgGSII$, a simple manipulation reveals that \eqref{eq:1-ptt} is equivalent to
\begin{equation}
  \label{eq:2-ptt}
  \left(
    \exponential{ p \Tr \lcgGSII } - 1
  \right)
  \Tr
  \left(
    \lcgGSII - 2 \identity + \inverse{\lcgGSII}
  \right)
  +
  \Tr
  \left(
    \inverse{\lcgGSII} + \ln \lcgGSII - \identity
  \right)
  \ge
  0.
\end{equation}
Since $p > 0$, and $\Tr \lcgGSII > 0$, the factor $\exponential{ p \Tr \lcgGSII } - 1$ is positive. Moreover, both trace terms in \eqref{eq:2-ptt} have been already shown to be nonnegative, see \eqref{eq:f-assumption-1-oldroyd-b} and \eqref{eq:3-oldroyd-b}. Requirement \eqref{eq:fenergy-1-f-stability-inequality-appendix} is thus fulfilled.

%%% Local Variables:
%%% mode: latex
%%% TeX-master: "../heat-stability-viscoelastic"
%%% End:

\vfill
\end{document}